\documentclass[11pt]{amsart}
\usepackage{amsopn,amsmath,amsthm,amssymb}
\usepackage{latexsym}
\usepackage[hypertex]{hyperref}

\newcommand{\nc}{\newcommand}

\nc{\vg}{\mathfrak{v} } \nc{\wg}{\mathfrak{w} }
\nc{\zg}{\mathfrak{z} } \nc{\ngo}{\mathfrak{n} }
\nc{\kg}{\mathfrak{k} } \nc{\mg}{\mathfrak{m} }
\nc{\bg}{\mathfrak{b} } \nc{\ggo}{\mathfrak{g} }
\nc{\rg}{\mathfrak{r} } \nc{\ggob}{\overline{\mathfrak{g}} }
\nc{\sog}{\mathfrak{so} } \nc{\sug}{\mathfrak{su} }
\nc{\spg}{\mathfrak{sp} } \nc{\slg}{\mathfrak{sl} }
\nc{\glg}{\mathfrak{gl} } \nc{\cg}{\mathfrak{c} }
\nc{\hg}{\mathfrak{h} } \nc{\tg}{\mathfrak{t} }
\nc{\ug}{\mathfrak{u} } \nc{\dg}{\mathfrak{d} }
\nc{\ag}{\mathfrak{a} } \nc{\osc}{\mathfrak{osc} }
\nc{\pg}{\mathfrak{p} } \nc{\sg}{\mathfrak{s} }
\nc{\pca}{\mathcal{P}} \nc{\nca}{\mathcal{N}}

\nc{\vp}{\varphi} \nc{\ddt}{\frac{{\rm d}}{{\rm d}t}}

\nc{\SO}{{\sf SO}} \nc{\Spe}{{\sf Sp}} \nc{\Sl}{{\sf Sl}}
\nc{\SU}{{\sf SU}} \nc{\Or}{{\sf O}} \nc{\U}{{\sf U}}
\nc{\Gl}{{\sf Gl}} \nc{\Se}{{\sf S}} \nc{\Cl}{{\sf Cl}}
\nc{\Spin}{{\sf Spin}} \nc{\Pin}{{\sf Pin}}

\nc{\RR}{{\mathbb R}} \nc{\HH}{{\mathbb H}} \nc{\CC}{{\mathbb C}}
\nc{\ZZ}{{\mathbb Z}} \nc{\FF}{{\mathbb F}} \nc{\NN}{{\mathbb N}}
\nc{\GG}{{\mathbb G}} \nc{\JJ}{{\mathbb J}} \nc{\II}{{\mathbb I}}
\nc{\KK}{{\mathbb K}} \nc{\DD}{{\mathbb D}}

\nc{\ad}{\operatorname{ad}} \nc{\Ad}{\operatorname{Ad}}
\nc{\rank}{\operatorname{rank}} \nc{\Irr}{\operatorname{Irr}}
\nc{\End}{\operatorname{End}} \nc{\Aut}{\operatorname{Aut}}
\nc{\Inn}{\operatorname{Inn}} \nc{\Der}{\operatorname{Der}}
\nc{\Ker}{\operatorname{Ker}} \nc{\Iso}{\operatorname{I}}
\nc{\Le}{\operatorname{L}} \nc{\tr}{\operatorname{tr}}
\nc{\dif}{\operatorname{d}\!} \nc{\sen}{\operatorname{sen}}
\nc{\modu}{\operatorname{mod}} \nc{\Ric}{\operatorname{R}}
\nc{\Sym}{\operatorname{Sym}} \nc{\sca}{\operatorname{sc}}
\nc{\scalar}{{\sf s}} \nc{\grad}{\operatorname{grad}}
\nc{\ricci}{\operatorname{r}} \nc{\riccin}{\operatorname{Ric}}
\nc{\Lie}{\operatorname{L}} \nc{\tang}{\operatorname{T}}
\newcommand\aff{{\mathfrak{aff}}}

\nc{\coad}{\operatorname{coad}} \nc{\lra}{\longrightarrow}
\nc{\e}{{\rm e}} \nc{\alt}{\raise1pt\hbox{$\bigwedge$}}

\theoremstyle{plain}
\newtheorem{thm}{Theorem}[section]
\newtheorem{prop}[thm]{Proposition}
\newtheorem{cor}[thm]{Corollary}
\newtheorem{lem}[thm]{Lemma}

\theoremstyle{definition}

\theoremstyle{remark}
\newtheorem*{rem}{Remark}
\newtheorem*{rems}{Remarks}
\newtheorem{exam}[thm]{Example}

\newcommand{\ri}{{\rm (i)}}
\newcommand{\rii}{{\rm (ii)}}
\newcommand{\riii}{{\rm (iii)}}

\title[Lie bialgebras of complex type]{Lie bialgebras of complex type and associated Poisson Lie groups}

\author{A. Andrada, M. L. Barberis} 
\address{CIEM, FaMAF, Universidad Nacional de C\'ordoba, Ciudad Universitaria, (5000) C\'or\-do\-ba, Argentina}
\email{andrada@mate.uncor.edu, barberis@mate.uncor.edu}
\author{G. Ovando}
\address{CIEM, FaMAF, Universidad Nacional de C\'ordoba, Ciudad Universitaria, (5000) C\'or\-do\-ba, Argentina - and - FCEIA, Universidad Nacional de Rosario, Av. Pellegrini 250, (2000) Rosario, Argentina}
\email{ovando@mate.uncor.edu}

\begin{document}

\begin{abstract}
In this work we study a particular class of Lie bialgebras arising
from Hermitian structures on Lie algebras such that the metric is ad-invariant. We
will refer to them as Lie bialgebras  of complex type. These give
rise to Poisson Lie groups $G$ whose corresponding duals $G^*$ are
complex Lie groups. We also prove that a Hermitian structure on
$\ggo$ with ad-invariant metric induces a structure of the same type
on the double Lie algebra ${\mathcal
D}\mathfrak{g}=\mathfrak{g}\oplus\mathfrak{g}^*$, with respect to
the canonical ad-invariant metric of neutral signature on ${\mathcal
D}\mathfrak{g}$. We show how to construct a $2n$-dimensional Lie
bialgebra of complex type starting with one of dimension $2(n-2), \;
n\geq 2$. This allows us to determine all solvable Lie algebras of
dimension $\leq 6$ admitting a Hermitian structure with ad-invariant
metric. We exhibit some examples in dimension $4$ and $6$, including
two one-parameter families, where we identify the Lie-Poisson
structures on the associated simply connected Lie groups, obtaining
also their symplectic foliations.
\end{abstract}

\maketitle

\section{Introduction}

The notion of Poisson Lie group was first  introduced by Drinfel'd
\cite{D}  and studied by Semenov-Tian-Shansky \cite{S-T-S} to
understand the Hamiltonian structure of the hidden symmetry group of
a completely integrable system. These groups appear in the theory of
Poisson-Lie T-dual sigma models \cite{Kl}. Lie bialgebras, the
infinitesimal counterpart of Poisson Lie groups (see \cite{D,W}),
are in a natural one-to-one correspondence with Lie algebra
structures on the vector space $\mathcal D \ggo = \ggo \oplus \ggo
^*$ with some compatibility conditions. $\mathcal D \ggo$ with this
Lie algebra structure is called the double of the Lie algebra
$\ggo$. The classification of Lie bialgebras for complex semisimple
$\ggo$ was carried out in \cite{BDr}. On the other hand, since
non-semisimple Lie algebras play an important role in physical
problems, there is more recent work on the classification problem of
low-dimensional Lie bialgebras (see, for instance, \cite{RHR} and
the references therein).

The purpose of this paper is to study a particular class of Lie
bialgebras arising from classical r-matrices of a special type and
to present some new examples of Poisson Lie groups in this setting.
We  obtain Poisson Lie groups such that their duals turn out to be
complex Lie groups. In fact, if $J$ is a complex structure on a Lie
algebra $\ggo$, then  a new bracket $[\cdot , \cdot ]_J$  can be
defined on $\ggo$ as follows:
\[ [x, y]_J = [Jx,y] + [x,Jy], \qquad  x, y \in \ggo,\]
and the integrability condition for $J$ (see \S\ref{cx-r}) implies
that $[\cdot , \cdot ]_J$ satisfies the Jacobi identity, that is,
$J$ is a classical r-matrix on $\ggo$ (see Lemma~\ref{complex}).
Moreover, $\ggo_J:=(\ggo, [\cdot, \cdot]_J)$ is a complex Lie
algebra. If $\ggo$ carries also an ad-invariant metric $g$ such
that $(J,g)$ is a Hermitian structure, then $\ggo$ admits a Lie
bialgebra structure. Such a Lie bialgebra will be called of {\em
complex type}. We prove that the Lie algebra structure on the dual
$\ggo^*$ is isomorphic to $\ggo_J$, and therefore $\ggo^*$ is a
complex Lie algebra. This implies that the simply connected
Poisson Lie group $G$ associated to the Lie bialgebra $\ggo$ has
the remarkable property that its dual Poisson Lie group $G^*$ is a
complex Lie group. Furthermore, we prove that the double Lie
algebra ${\mathcal D}\ggo$ admits a complex structure ${\mathcal
J}$ (induced by $J$), such that $({\mathcal J}, \langle\cdot\,
,\cdot\rangle)$ is a Hermitian structure on ${\mathcal D}\ggo$,
where $\langle\cdot\, ,\cdot\rangle$ is the canonical ad-invariant
metric of neutral signature on ${\mathcal D}\ggo$. It turns out
that when $\ggo$ is a Lie bialgebra of complex type, then so is
${\mathcal D}\ggo$.

Also, beginning with a Lie algebra $\ggo$ equipped with a Hermitian
structure such that the metric is ad-invariant, we can apply, under
certain conditions, the double extension method (see \cite{FS, Ka,
MR}) in order to construct a new Lie algebra of dimension
$\dim\ggo+4$ endowed with a similar structure. This construction
leads us to determine all solvable Lie algebras of dimension $\leq
6$ admitting a Hermitian structure with ad-invariant metric.

The paper is organized as follows. In \S2.1 we define Lie
bialgebras,  Poisson Lie groups and recall from \cite{AGMM} how to
obtain a Lie bialgebra from a classical $r$-matrix. Section~2.2 is
devoted to preliminary results on the complex Lie algebra $\ggo_J$
constructed from a complex structure $J$ on $\ggo$ and in \S2.3 we
define metric Lie algebras and discuss some properties of Hermitian
structures with ad-invariant metrics.

In Section 3 we prove that every Hermitian structure with
ad-invariant metric gives rise to a Lie bialgebra of complex type.
In this case, we show that the double Lie algebra $\mathcal D \ggo$
also has such a structure. We obtain as a corollary  that a complex
structure on $\ggo$ induces a Lie bialgebra structure of complex
type on the double $\mathcal D (T^*\ggo)$ of the cotangent Lie
algebra $T^*\ggo$. We exhibit many examples starting with solvable
Lie algebras. In case $\ggo$ is compact semisimple, we show that our
construction yields a Lie bialgebra which is not equivalent to the
standard one obtained by \cite{LW}.

Section 4 is devoted to the study of Lie bialgebras of complex type
obtained by a double extension process. By applying this
construction we can prove that every solvable metric Lie algebra of
dimension at most $6$ with  metric of signature $(2r,2s)$ admits a
Lie bialgebra structure of complex type.

In the last section we study the simply connected Poisson Lie
groups associated to some Lie bialgebras of complex type in
dimension $4$ and $6$, obtaining also their symplectic foliations.
We get a linearizable Poisson structure on $\RR^4$ and we exhibit
two one-parameter families $\{ \Pi_1^{\lambda}\}$, $\{
\Pi_2^{\lambda}\}$ of Lie-Poisson structures on $\RR^6$ ($\lambda
>0$) such that $\Pi_1^{\lambda}$ is not equivalent to
$\Pi_2^{\lambda}$ for each $\lambda >0$, but as $\lambda \to 0$
these families converge to the same Poisson structure $\Pi^0$. When
$\lambda \to \infty$ we have an analogous situation, with both
families converging to the same Poisson structure $\Pi^{\infty}$. We
point out that the structures $\Pi^0$ and $\Pi^{\infty}$ on $\RR^6$
are not equivalent. On the other hand, there is a unique nilpotent
Lie algebra of dimension 6 admitting Lie bialgebra structures of
complex type, giving rise to Lie-Poisson structures on $\RR^6$ which
cannot be obtained as deformations of the previous ones.

\

\section{Preliminaries and basic results}

\subsection{Lie bialgebras and Poisson Lie groups}
Let $\ggo$ be a real Lie algebra and $\delta : \ggo \to \alt ^2
\ggo $ a $1$-cocycle with respect to the adjoint representation.
The pair $(\ggo , \delta )$ is called a {\em Lie bialgebra} if
$\delta^*:\alt ^2 \ggo^*\to \ggo ^* $ induces a Lie algebra
structure on $\ggo ^*$ (see \cite{D}).

A \emph{multiplicative} Poisson structure on a Lie group $G$ is a
smooth section $\Pi$ of $\Lambda ^2 TG$ such that the following
conditions are satisfied:
\begin{enumerate}
  \item $\{f,g\} := \langle df \wedge dg , \Pi\rangle$ defines a Lie
  bracket on $C^{\infty}(G, \RR)$;
  \item the multiplication $m:G\times G \to G$ is a Poisson mapping,
  that is, the pull back mapping $m^*: C^{\infty}(G, \RR) \to C^{\infty}(G \times G , \RR)$
  is a homomorphism for the Poisson brackets\footnote{We recall that the Poisson bracket on
  $C^{\infty}(G \times G,\RR)$ is defined by
\[ \{f,g\}(x,y) =\{f(x,\cdot ), g(x,\cdot )\}(y)+ \{f(\cdot,y), g(\cdot ,y)\}(x),\qquad x,y \in G.
\]}.
\end{enumerate}
A Lie group equipped with a multiplicative Poisson structure is
called a \emph{Poisson Lie group}.

If $G$ is a connected, simply connected Lie group with Lie algebra
$\ggo$, there is a one-to-one correspondence between
multiplicative Poisson structures on $G$ and Lie bialgebra
structures $\delta$ on $\ggo$ (see \cite{Dr}). The correspondence
is established as follows:
\begin{equation}\label{corresp}
    (d\Pi )_e = \delta .
\end{equation}

We will say that a linear operator $r \in \End(\ggo)$ is  a {\it
classical r-matrix} if the $\ggo$-valued skew-symmetric bilinear
form on $\ggo$ given by
\[ [x, y]_r = [rx,y] + [x,ry]\] is a Lie bracket, that is, it
satisfies the Jacobi identity. $\ggo$ equipped with the Lie
bracket $[\cdot , \cdot ]_r$ will be denoted by $\ggo _r$. The
next result is a consequence of Corollary~2.9 in \cite{AGMM} (see
also \cite{STS}):

\begin{prop} \label{r-bialg}
If $r$ is a classical r-matrix and $\ggo$ admits a non-degenerate
symmetric ad-invariant bilinear form such that $r$ is
skew-symmetric, then $r$ gives rise to a Lie bialgebra $(\ggo ,
\delta _r )$.
\end{prop}

The 1-cocycle $\delta _r$ above is defined as follows. Let $g$ be a
non-degenerate symmetric ad-invariant bilinear form on $\ggo$ such
that $r$ is skew-symmetric. We consider the metric $g$ on $\ggo$ as
a linear isomorphism $g:\ggo\lra\ggo^*,\;g(x)(y)=g(x,y)$. Via this
isomorphism, the Lie bracket $[\cdot,\cdot]_r$ on $\ggo_r$ induces a
Lie bracket $[\cdot,\cdot]^*$ on $(\ggo_r)^*=\ggo^*$, given as
follows:
\[ [\alpha,\beta]^*=g([g^{-1}\alpha,g^{-1}\beta]_r), \qquad
\alpha,\beta\in\ggo^*.
\]  We denote by $\delta_r:
\ggo \to \alt ^2 \ggo $ the dual of $[\cdot,\cdot]^*$. If we look
upon $R:=r \circ g^{-1}$ as an element in $\alt^2 \ggo $, we
observe that $(\ggo , \delta _r )$ is an exact Lie bialgebra, since
\[ \delta _r (x) = \ad (x) R\qquad \text{ for all } x
\in \ggo,\] that is, $\delta _r$ is a coboundary.

For a Lie bialgebra coming from a classical r-matrix $r$, the
corresponding Poisson structure on $G$ given by \eqref{corresp} is
the bivector field
$\stackrel{\leftarrow}{R}-\stackrel{\rightarrow}{R}$, where
$\stackrel{\leftarrow}{R}$ (resp. $\stackrel{\rightarrow}{R}$) is
the left (resp. right) invariant bivector field on $G$ whose value
at the identity $e$ of $G$ is $R=r\circ g^{-1}$ (see
\cite{AGMM,CP}).

In \S\ref{bi-herm} we will apply the above proposition to the
particular case when the r-matrix is a complex structure on the Lie
algebra.

\begin{rem} Given a bivector $R \in \alt^2 \ggo$, the coboundary
$\delta_R$ defined by $\delta _R (x) = \ad (x) R$ for $x \in \ggo$
induces a Lie bialgebra structure on $\ggo$ if and only if the
Schouten bracket\footnote{For the definition and properties of the
Schouten bracket see, for instance, \cite{K}.} $[R,R]$ is
$\ad(\ggo)$-invariant:
\begin{equation*} [R,R]\in \left(\alt^3 \ggo\right)^{\ggo}.\end{equation*}
When $R$ satisfies this condition we say that $R$ is a solution of
the {\it modified Yang-Baxter equation} (MYBE). In particular,
this is trivially satisfied when $[R,R]=0$, which is known as the
{\it classical Yang-Baxter equation} (CYBE).

If $\ggo$ admits a metric as above and $r=R\circ g \in \End(\ggo)$,
then the MYBE and the CYBE in terms of $r$ are equivalent,
respectively, to:
\begin{gather*}
[x, B_r(y,z)]+[y,B_r(z,x)]+[z,B_r(x,y)]=0, \qquad \quad x, y, z
\in \ggo,\\ B_r(x,y) =0, \qquad \quad x, y \in \ggo ,
\end{gather*}
where
\[ B_r(x,y):= [rx,ry]-r\left([x,y]_r\right) .\]
\end{rem}

\medskip

\subsection{Complex structures as $r$-matrices and associated Lie algebras}\label{cx-r}
A complex structure on a real Lie
algebra $\ggo$ is an endomorphism $J$ of $\ggo$ satisfying $J^2=-$id
and the integrability condition $N_J \equiv 0$, where
\begin{equation}\label{integ} N_J(x,y)=[Jx,Jy]-[x,y]-J[Jx,y]-J[x,Jy] \qquad \text{ for }
x, y \in \ggo.\end{equation} If $G$ is a Lie group with Lie
algebra $\ggo$, by left translating the endomorphism $J$ we obtain
a  complex manifold $(G,J)$ such that left translations are
holomorphic maps. We point out that $(G,J)$ is not necessarily a
complex Lie group since right translations are not in general
holomorphic.

Let $\ggo^{\CC}$ denote the complexification of $\ggo$, then we
have a splitting
\[\ggo^{\CC}=\ggo^{1,0}\oplus \ggo^{0,1},\]
where $\ggo^{1,0}$ (resp. $\ggo^{0,1}$) is the eigenspace of $J$ of
eigenvalue $i$ (resp. $-i$). It follows that
\[ \ggo^{1,0}=\{ x-iJx : x \in \ggo \}, \qquad \ggo^{0,1}=  \{ x+iJx : x \in \ggo \},   \]
hence $\ggo^{0,1}=\sigma (\ggo^{1,0})$, where $\sigma : \ggo
^{\CC} \to \ggo ^{\CC}$ is conjugation with respect to the real
form $\ggo$, that is, $\sigma (x + i y) = x - i y$, $x, y \in
\ggo$. The integrability of $J$ is equivalent to the fact that
both, $\ggo^{1,0}$ and $\ggo^{0,1}$, are complex Lie subalgebras
of $\ggo^{\CC}$. When $\ggo^{1,0}$ is abelian $J$ is called an
{\em abelian} complex structure  (see \cite{BD2}) and it follows
that $\ggo^{0,1}$ is also abelian. This is equivalent to the
condition $[Jx,Jy]=[x,y], \; x, y \in \ggo.$

When $\ggo^{1,0}$ is an ideal of $\ggo^{\CC}$, then $\ggo^{0,1}$
is also an ideal and $(\ggo,J)$ is called a complex Lie algebra.
In this case, $\ad(x)$, $x \in \ggo$, are complex linear maps:
\begin{equation}\label{cxlin}
\ad (x)\circ J =J\circ \ad (x).  \end{equation} If $G$ is a
connected Lie group with Lie algebra $\ggo$ and $(G,J)$ is the
complex manifold obtained by left translating $J$, the above
condition is equivalent to $(G,J)$ being a complex Lie group, that
is, both, right and left translations on $G$ are holomorphic maps.

We will say that two complex structures $J_1, \, J_2$ on $\ggo$
are equivalent if there exists a Lie algebra automorphism
$\varphi:\ggo\lra\ggo$ such that $\varphi \circ J_1=J_2 \circ
\varphi$.

As a consequence of the integrability condition of a complex
structure $J$ we obtain the  following lemma (see also
\cite{Ba,LS}).

\begin{lem} \label{complex}
If $J$ is a complex structure on $\ggo$ then $J$ is a classical
r-matrix. Moreover, $(\ggo_J , J)$ is a complex Lie algebra.
\end{lem}

\begin{proof} Using \eqref{integ} we compute
\begin{eqnarray*} J[x,y]_J &=& J([Jx,y]+[x,Jy])=-[x,y]+[Jx,Jy], \\
 {[Jx,y]_J} &=& [J(Jx),y]+[Jx,Jy]=-[x,y]+[Jx,Jy],\end{eqnarray*}
therefore
\begin{equation} \label{J}
 J[x,y]_J=[Jx,y]_J.
\end{equation}
Setting $\ad_J(x)=[x,\cdot ]_J$, \eqref{J} amounts to
$\ad_J(x)\circ J = J\circ \ad_J(x)$, for all $x \in \ggo_J$, in
other words, $(\ggo_J,J)$ is a complex Lie algebra and the lemma
follows.
\end{proof}

\begin{lem}\label{lem1}
\begin{enumerate} \item [ ]
\item $\ggo_J$ is an abelian Lie algebra if and only if $J$ is an
abelian complex structure. \item  $(\ggo,  J)$ is a complex Lie
algebra if and only if it is isomorphic to $(\ggo _J, J)$. \item
$(\ggo, J)$ is a complex Lie algebra if and only if $\ad_J(x)= 2 J
\circ \ad(x)$ for any $x \in \ggo$.
\end{enumerate}\end{lem}

\begin{proof}
Assertion (1) is straightforward from the definition of
$[\cdot,\cdot ]_J$.

To prove (2), assume first that $(\ggo,J)$ is a complex Lie
algebra. Then \[ [x,y]_J =2J[x,y], \qquad x,y \in \ggo.\] It
follows that $\varphi=-\dfrac12 J$ is a Lie algebra isomorphism
from $\ggo$ onto $\ggo_J$ that commutes with $J$.

For the converse, let $\varphi$ be a Lie algebra isomorphism
$\varphi : \ggo_J \to \ggo$ such that $\varphi \circ J = J \circ
\varphi$. Then
\begin{eqnarray*}
\ad (x) \circ J \;(y)&=&[x, Jy]= \varphi \left( [\varphi ^{-1}(x),
\varphi ^{-1}(Jy)]_J  \right)= \varphi \left([\varphi ^{-1}(x), J
\circ \varphi ^{-1}(y)]_J \right)\\ &=& \varphi \circ J \left(
[\varphi ^{-1}(x), \varphi ^{-1}(y)]_J \right) = J \circ
\varphi  \left( [\varphi ^{-1}(x), \varphi ^{-1}(y)]_J \right) \\
&=& J \left([x,y]\right) = J\circ \text{ad} (x)\; (y),
\end{eqnarray*}
where the fourth equality follows from Lemma~\ref{complex}.
Therefore, $(\ggo,J)$ is a complex Lie algebra.

The only if part of (3) is straightforward. To prove the converse,
assume that $\ad_J(x)= 2 J \circ \ad(x)$ for any $x \in \ggo$.
Since $J$ is integrable, we have
\[ [x,y]_J=2 J[x,y]=2\left([Jx,y]+[x,Jy]+J[Jx ,Jy] \right),
\]
hence, \[ 2J[Jx,Jy]= -[Jx,y] -[x,Jy]=-2J[x,y], \] which implies
$[Jx,Jy]=-[x,y]$. In particular, $[Jx,y]=[x,Jy]$, and the
integrability of $J$ gives $J[x,y]=[Jx,y]$, that is, $(\ggo,
 J)$ is a complex Lie algebra.
\end{proof}

The next proposition gives a complex Lie algebra isomorphism between
$(\ggo_J, J)$ and $\ggo ^{1,0}$ (see also \cite{LS}).

\begin{prop}\label{prop1}
If $J$ is a complex structure on $\ggo$, then $(\ggo_J, J)$ and
$\ggo ^{1,0}$ are isomorphic as complex Lie algebras.
\end{prop}

\begin{proof}
Consider $\ggo$ as a complex vector space with multiplication by
$i$ given by the endomorphism $J$. Let $\varphi: \ggo _J \to
\ggo^{1,0}$ be the  complex linear map defined  as follows
\[ \varphi (x)= i(x -iJx)=(J +i \, \text{id})(x), \qquad x\in \ggo _J,  \]
where id is the identity operator. We need to check that $\varphi$
is a Lie algebra homomorphism, where the Lie bracket on $\ggo^{1,0}$
is the complex bilinear extension of $[\cdot , \cdot ]$. In fact,
\begin{eqnarray*} {\varphi ([x , y ] _J}) &=& \varphi([Jx,y]+[x ,Jy])
= J([Jx,y]+[x ,Jy]) +i([Jx,y]+[x ,Jy]) \\
&=& {[Jx , Jy]-[x,y]+  i([Jx,y]+[x ,Jy])}, \\
 {[\varphi (x) ,\varphi (y) ] }&=& {[Jx +ix, J +iy]}
 ={[Jx , Jy]-[x,y] +i([Jx,y]+[x ,Jy])},
\end{eqnarray*}
where the third equality follows from \eqref{integ}. Hence,
$\varphi$ establishes the desired Lie algebra isomorphism.
\end{proof}

\begin{rems}\

\smallskip

\noindent $\ri$  Observe that in the above proposition $J-i\,
\text{id}$ gives an isomorphism between $\ggo_J$ and $\ggo^{0,1}$.

\smallskip

\noindent $\rii$ We point out  that if  $\ggo$ is a solvable (resp.
nilpotent) Lie algebra with a complex structure $J$, then $\ggo_J$
is solvable (resp. nilpotent). The converse  is not true, that is,
we can obtain a solvable Lie algebra $\ggo_J$ starting from a non
solvable Lie algebra $\ggo$. We recall from \cite{SD} the definition
of a regular complex structure on a reductive Lie algebra $\ggo$. A
complex structure $J$ on $\ggo$ is said to be {\em regular} if there
exists a $\sigma$-stable Cartan subalgebra $\hg$ of $\ggo^{\CC}$
such that $[\hg, \ggo ^{1,0}] \subset \ggo^{1,0}$, where $\sigma$ is
conjugation in $\ggo^{\CC}$ with respect to $\ggo$. If $\ggo$ is a
product of an abelian ideal and a semisimple ideal $\sg$, it was
shown in \cite{SD} that when the simple factors in $\sg$ are compact
or belong to the following list:
\[ \mathfrak{sl}(2,\RR), \;\; \mathfrak{sp}(2n,\RR),\;\; \mathfrak{so}^*(2n),
\;\;  \mathfrak e _{6(-14)}, \;\; \mathfrak e _{6(2)}, \;\;
\text{any real form of: }  \mathfrak e _7, \;\;\mathfrak e_8, \;\;
\mathfrak f_4, \;\; \mathfrak g_2 , \] then any complex structure
$J$ on $\ggo$ is regular and satisfies $\ug \subset \ggo^{1,0}
\subset \mathfrak b$ for some Borel subalgebra $\mathfrak b$  of
$\ggo^{\CC}$ with unipotent radical $\ug$. Therefore, $\ggo_J\cong
\ggo^{1,0}$ is always solvable.

On the other hand, if $\ggo$ is a a real noncompact semisimple Lie
algebra with a compact Cartan subalgebra or a real split
semisimple Lie algebra, it follows from results of \cite{LQ} that
for the canonical Koszul operator $J$ on $\ggo$, $\ggo_J$ is also
solvable.
\end{rems}

The next result is a straightforward consequence of
Proposition~\ref{prop1}.

\begin{cor}
If $J_1$ and $J_2$ are equivalent complex structures on $\ggo$,
then $\ggo _{J_1}$ is isomorphic to $\ggo _{J_2}$.
\end{cor}

\medskip

\subsection{Metric Lie algebras and Hermitian structures }
A {\em metric} on $\ggo$ is a non degenerate symmetric bilinear
form. A Hermitian structure $(J,g)$ on $\ggo$ is a pair of a complex
structure $J$ and a metric $g$ on $\ggo$ such that $J$ is
skew-symmetric. We will say that two Lie algebras with Hermitian
structures $(\ggo_1,J_1,g_1)$ and $(\ggo_2,J_2,g_2)$ are {\em
equivalent} if there exists a Lie algebra isomorphism
$\varphi:\ggo_1\lra\ggo_2$ which is also an isometry between $g_1$
and $g_2$ satisfying $\varphi \circ J_1=J_2 \circ \varphi$.

A {\em metric Lie algebra} is a pair $(\ggo, g)$ of a Lie algebra
$\ggo$ equipped with an ad-invariant metric $g$, that is, the
following condition holds for all $x,y,z\in\ggo$:
\[ g([x,y],z)=g(x,[y,z]).\]
Equivalently, the endomorphisms $\ad(x)$ are skew-symmetric for
all $x\in\ggo$. Note that in this case $\ggo$ is unimodular (i.e.,
$\tr(\ad(x))=0$ for all $x\in\ggo$). The Killing form on a
semisimple Lie algebra is an example of an ad-invariant metric.
However, for non semisimple Lie algebras, there are obstructions
for the existence of ad-invariant metrics: for instance,  if such
a metric exists, then $\zg^{\perp}=[\ggo,\ggo]$, where $\zg$
denotes the centre of $\ggo$. In particular, we have that
$\dim\zg+\dim[\ggo,\ggo]=\dim\ggo$, so that a solvable Lie algebra
with trivial centre cannot admit an ad-invariant metric. The
metric Lie algebra $(\ggo ,g)$ is called {\em indecomposable} if
every proper ideal of $\ggo$ is degenerate.

We will be studying Hermitian structures $(J,g)$ on $\ggo$ where
$g$ is ad-invariant. First, we prove some results which impose
restrictions on the complex structure or the metric. The next
lemma gives a necessary condition on the signature of the metric.
It can be shown by induction on $\dim V$ and its proof is omitted.

\begin{lem}\label{sign}
Let $V$ be a real vector space with a Hermitian structure $(J,g)$.
Then the signature of $g$ is of the form $(2r,2s),\,2r+2s=\dim V$.
\end{lem}

\smallskip

\begin{prop}
Let $\ggo$ be a Lie algebra endowed with a Hermitian structure $(J,
g)$ such that the metric  $g$ is ad-invariant. If any of the
conditions $\ri-\riii$ below holds, then $\ggo$ is abelian.
\begin{enumerate}
\item[\ri] The K\"ahler form $\omega$ on $\ggo$ defined by
$\omega(x,y)=g(Jx,y)$ is closed. \item[\rii] $J$ is abelian.
\item[\riii] $J \circ \ad(x)=\ad(x) \circ J$ for all $x\in\ggo$.
\end{enumerate}
\end{prop}

\begin{proof}
$\ri$ Suppose that the $2$-form $\omega$ is closed, so that
\[ \omega([x,y],z)+\omega([y,z],x)+\omega([z,x],y)=0 \qquad \text{for all } x,y,z\in\ggo.\]
Since $g$ is ad-invariant and $J$ skew-symmetric, from the
definition of $\omega$ we obtain
\[ \omega([x,y],z)=g(J[x,y],z),\,\omega([y,z],x)=g(z,[y,Jx]),\, \omega([z,x],y)=-g(z,[x,Jy]),\]
so that
\[ g(J[x,y]-[Jx,y]-[x,Jy],z)=g(J[Jx,Jy],z)=0 \]
for all $x,y,z\in\ggo$, where we have used the integrability of
$J$. Since $g$ is non degenerate, we have that $[Jx,Jy]=0$ for all
$x,y\in\ggo$, and hence $\ggo$ is abelian.

\smallskip

$\rii$ Let us assume now that $J$ is abelian. Then
\[ g([Jx,y],z)=-g(y,[Jx,z])\stackrel{\dagger}{=}g(y,[x,Jz])=g(J[x,y],z),\]
where in equality $(\dagger)$ we have used the fact that $J$ is
abelian. Since $g$ is non degenerate, we have that $J[x,y]=[Jx,y]$
for all $x,y\in\ggo$, and from this we obtain easily that
$[Jx,y]=[x,Jy]$. However, $[Jx,y]=-[x,Jy]$ as $J$ is abelian, and
hence $[Jx,y]=0$ for all $x,y$ and $\ggo$ is abelian.

\smallskip

$\riii$ Suppose now that $J\circ \ad(x)=\ad(x)\circ J$, or
equivalently, $J[x,y]=[Jx,y]$ for all $x,y\in\ggo$. Then we have
\[ g([Jx,y],z)=-g(y,[Jx,z])=-g(y,J[x,z])=-g([x,Jy],z),\]
so that $[Jx,y]=-[x,Jy]$ and therefore $J$ is abelian. From $\rii$,
$\ggo$ is abelian.
\end{proof}

\

\section{Lie bialgebras and Hermitian structures}\label{bi-herm}

Let $\ggo$ be a Lie algebra with a Hermitian structure $(J,g)$. We
consider the metric $g$ on $\ggo$ as a linear isomorphism
$g:\ggo\lra\ggo^*,\;g(x)(y)=g(x,y)$. Via this isomorphism, the Lie
bracket $[\cdot,\cdot]_J$ on $\ggo_J$ induces a Lie bracket
$[\cdot,\cdot]^*$ on $(\ggo_J)^*=\ggo^*$, given as follows:
\[ [\alpha,\beta]^*=g([g^{-1}\alpha,g^{-1}\beta]_J), \qquad \alpha,\beta\in\ggo^*
\] and $J$ induces an endomorphism $J^*:\ggo ^* \to \ggo ^*$,
\begin{equation} \label{J*} J^*\alpha=-\alpha\circ J, \qquad \alpha \in \ggo ^* .
\end{equation}
We denote by $\delta_J: \ggo \to \ggo \wedge \ggo$ the dual of
$[\cdot,\cdot]^*$.

\begin{thm} \label{bialg}
Let $(J,g)$ be a Hermitian structure on $\ggo$ such that $g$ is
ad-invariant. Then $(\ggo,\delta _J)$ is a Lie bialgebra such that
$(\ggo ^* , [\cdot , \cdot ]^*, J^*)$ is a complex Lie algebra
isomorphic to $(\ggo_J ,[\cdot, \cdot]_J, J)$. In particular, if
$G$ and $G^*$ are the corresponding simply connected Poisson Lie
groups, then $G^*$ is a complex Lie group.
\end{thm}

\begin{proof} The fact that $(\ggo,\delta _J)$ is a Lie bialgebra
is a consequence of Proposition~\ref{r-bialg}. It follows from
Proposition 2.8 in \cite{AGMM} that
\[ (\ggo^*,[\cdot ,\cdot ]^*) \cong (\ggo_J ,[\cdot, \cdot]_J).\]
Moreover,  since $ J^*\circ g=g \circ J$, Lemma~\ref{complex}
implies that $(\ggo^*,[\cdot ,\cdot ]^*,J^*) \cong (\ggo_J ,[\cdot,
\cdot]_J,J)$ as complex Lie algebras.
\end{proof}

A Lie bialgebra of the form $(\ggo,\delta _J)$ as above will be
called a {\it Lie bialgebra of complex type} and a Poisson Lie group
corresponding to it, a {\it Poisson Lie group of complex type}.

\begin{exam}\label{cot}
The {\em cotangent Lie algebra} $T^*\hg$ of a Lie algebra $\hg$ is
defined as the semidirect product $T^*\hg := \hg \ltimes_{\coad}
\hg^*$. If $H$ denotes the simply connected Lie group with Lie
algebra $\hg$, the cotangent bundle $T^*H$ has a natural Lie group
structure with corresponding Lie algebra $T^*\hg$. There is a
canonical ad-invariant metric of neutral signature $\langle
\cdot\, , \cdot \rangle$ on $T^*\hg$, defined as in \eqref{bil}
below. According to Theorem~\ref{bialg}, every skew-symmetric
complex structure  on $T^*\hg$ determines a multiplicative Poisson
structure on $T^*H$ such that $\left(T^*H\right)^*$ is a complex
Lie group. It should be pointed out that  Hermitian structures
$(J, \langle \cdot\, , \cdot \rangle)$ on $T^*\hg$ are in
one-to-one correspondence with left invariant generalized complex
structures on $H$, as shown in \cite{ABDF}.
\end{exam}

\

It is well known \cite{D} that given a Lie bialgebra $(\ggo,
\delta)$, the vector space ${\mathcal D}\ggo=\ggo\oplus\ggo^*$
becomes a Lie algebra with the bracket:
\[ [(x,\alpha),(y,\beta)]=\left([x,y]+\coad^*(\alpha)y-\coad^*(\beta)x,
[\alpha,\beta]^*+\coad(x)\beta-\coad(y)\alpha\right),\] for
$x,y\in\ggo,\,\alpha,\beta\in\ggo^*$, where $\coad$ stands for the
coadjoint representation of $\ggo$ on $\ggo^*$, and $\coad^*$
stands for the coadjoint representation of $\ggo^*$ on $\ggo$.
Note that both $\ggo$ and $\ggo^*$ are subalgebras of ${\mathcal
D}\ggo$ with this Lie bracket. There is a canonical ad-invariant
metric $\langle \cdot \, , \cdot \rangle$ of neutral signature on
${\mathcal D}\ggo$ given by
\begin{equation} \label{bil}\langle
(x,\alpha),(y,\beta)\rangle=\alpha(y)+\beta(x)\end{equation} for
$x,y\in\ggo,\;\alpha,\beta\in\ggo^*$. Given an almost complex
structure $J$ on $\ggo$ we define an almost complex structure
${\mathcal J}$ on ${\mathcal D}\ggo$  by
\begin{equation}\label{cx-str} {\mathcal J}(x,\alpha)=(Jx,J^*\alpha),\end{equation}
where $J^*$ is defined as in \eqref{J*}. Note that $\mathcal J$ is
skew-symmetric with respect to the metric \eqref{bil} and that
both subalgebras, $\ggo$ and $\ggo^*$, are ${\mathcal
J}$-invariant.

\smallskip

When the bialgebra structure on $\ggo$ is trivial, that is,
$\delta=0$, the Lie algebra ${\mathcal D}\ggo$ is the cotangent
Lie algebra $T^*\ggo$. It was proved in \cite{BD3} that if $J$ is
a complex structure  on $\ggo$ and we define $\mathcal J$ as in
\eqref{cx-str}, then $\mathcal J$ is a complex structure on the
cotangent algebra $T^*\ggo$. Moreover, $(\mathcal J ,\langle \cdot
\, , \cdot \rangle)$ is a Hermitian structure on $T^*\ggo$, where
$\langle \cdot \, , \cdot \rangle$ is the canonical ad-invariant
neutral  metric. More generally, we have:

\begin{prop} \label{double}
Let $\ggo$ be a Lie algebra with a Hermitian structure $(J,g)$ such
that $g$ is ad-invariant. Then the almost complex structure
$\mathcal J$ on the double Lie algebra ${\mathcal
D}\ggo=\ggo\oplus\ggo^*$ is integrable, with $\mathcal J$ defined as
in \eqref{cx-str}.
\end{prop}

\begin{proof}
Since $J$ is integrable on $\ggo$, $J^*$ is integrable on $\ggo^*$
(see Theorem~\ref{bialg}) and both are subalgebras of ${\mathcal
D}\ggo$, we only have to check that $N_{\mathcal
J}((x,0),(0,\beta))=0$ for all $x\in\ggo,\,\beta\in\ggo^*$. From the
definitions we arrive at
\begin{align*}
[{\mathcal J}(x,0),{\mathcal J}(0,\beta)] &
=(-\coad^*(J^*\beta)(Jx),\,\coad(Jx)(J^*\beta))\\
[(x,0),(0,\beta)] &
=(-\coad^*(\beta)x,\,\coad(x)\beta),\\
{\mathcal J}[{\mathcal J}(x,0),(0,\beta)] & =(-J\coad^*(\beta)(Jx),\,
J^*\coad(Jx)\beta),\\
{\mathcal J}[(x,0),{\mathcal J}(0,\beta)] &
=(-J\coad^*(J^*\beta)x,\,J^*\coad(x)J^*\beta).
\end{align*}
Let us begin with the first coordinate; then, for $\psi\in\ggo^*$,
we compute
\begin{align*}
(\coad^*(J^*\beta)(x))\psi & = (-x\circ \ad^*(J^*\beta))\psi = -x([J^*\beta,\psi]^*) \\
& = -x(J^*[\beta,\psi]^*) \quad\text{(since }(\ggo^*,J^*) \text{ is a complex Lie algebra)}\\
              & = -J^*[\beta,\psi]^*(x) = [\beta,\psi]^*(Jx)\\
& = (Jx)[\beta,\psi]^* = -(\coad^*(\beta)(Jx))\psi,
\end{align*}
and from this equality one easily verifies that the first
coordinate of $N_{\mathcal J}((x,0),(0,\beta))$ vanishes. Let us
compute now its second coordinate. For $v\in\ggo$, we have
\begin{align*}
(\coad(Jx)(J^*\beta))v & = -(J^*\beta)[Jx,v]  = \beta(J[Jx,v])\\
                 & = \beta([Jx,Jv]-[x,v]-J[x,Jv]) \quad\text{(using the integrability of $J$)}\\
                 & = J^*(\coad(Jx)\beta)v+(\coad(x)\beta)v+(J^*\coad(x)J^*\beta)v,
\end{align*}
and from this we see immediately that the second coordinate of
$N_{\mathcal J}((x,0),(0,\beta))$ also vanishes. Therefore,
$N_{\mathcal J}\equiv 0$.
\end{proof}

Proposition \ref{double} implies that $(\mathcal J , \langle \cdot
\, , \cdot \rangle)$ is a Hermitian structure on ${\mathcal
D}\ggo$, therefore, Theorem~\ref{bialg} yields:

\begin{cor}\label{bi-doble}  \begin{enumerate}\item[]
\item Under the same hypotheses of the previous proposition,
$({\mathcal D}\ggo , \delta_{\mathcal J})$ is a Lie bialgebra of
complex type. \item If $J$ is a complex structure on $\ggo$ and
$\langle \cdot \, , \cdot \rangle$ is the canonical ad-invariant
neutral metric on $T^*\ggo$, then $(\mathcal J ,\langle \cdot\, ,
\cdot \rangle)$ is a Hermitian structure on $T^*\ggo$. In
particular, $\mathcal D \left(T^*\ggo \right)$ is a Lie bialgebra
of complex type.
\end{enumerate}
\end{cor}

\begin{proof}
The first assertion is a consequence of Proposition~\ref{double}.
For (2), use the fact that $\mathcal J$ defined as in
\eqref{cx-str} is a complex structure on $T^*\ggo$ (see
\cite{BD3}) and therefore $(\mathcal J ,\langle \cdot\, , \cdot
\rangle)$ is a Hermitian structure on $T^*\ggo$. The second
assertion now follows by applying (1).
\end{proof}

\medskip

\begin{exam}\label{l2}
Let $\aff(\RR)$ be the non abelian real 2-dimensional Lie algebra
and let  $\ggo:=L_2(1,1)$ denote its cotangent Lie algebra (see
Example \ref{cot}), which has a basis $\{e_1,\ldots,e_4\}$ such that
the non vanishing Lie bracket relations are
\begin{equation}\label{cor-l2} [e_1,e_2]=e_2,\qquad [e_1,e_3]=-e_3,\qquad [e_2,e_3]=e_4.\end{equation}
Let $(J,g)$ be the Hermitian structure on $\ggo$ with ad-invariant
$g$ given by
\begin{gather*} Je_1=e_2,\qquad Je_3=-e_4,\\
g(e_1,e_4)=g(e_2,e_3)=1.
\end{gather*}

The Lie bracket $[\cdot,\cdot]_J$ is given by
\begin{equation}\label{c1} [e_1,e_3]_J=e_4,\qquad [e_1,e_4]_J=-e_3,\qquad [e_2,e_3]_J=e_3,
\qquad [e_2,e_4]_J=e_4, \end{equation} hence, $\ggo_J$ is
isomorphic, as a complex Lie algebra, to the complex analogue
$\aff(\CC)$ of $\aff(\RR)$. If $\{e^1,\ldots,e^4\}$ denotes the
dual basis of $\ggo^*$, from \eqref{c1} and the definition of $g$,
we have that the only non trivial brackets on $\ggo^*$ are \[
[e^1,e^3]^*=-e^1,\qquad [e^1,e^4]^*=e^2,\qquad [e^2,e^3]^*=-e^2,
\qquad [e^2,e^4]^*=-e^1.\] To determine the Lie algebra structure
on ${\mathcal D}\ggo$, it remains to compute the brackets of the
form $[e_i,e^j]=(-\coad^*(e^j)e_i,\,\coad(e_i)e^j)$, which is done
below:
\[ \begin{array}{lll}
[e_1,e^1]=-e_3,      & \;[e_2,e^1]=e_4,       &  \cr
[e_1,e^2]=-(e_4+e^2),  & \;[e_2,e^2]=-e_3+e^1,  &  \cr
[e_1,e^3]=e_1+e^3,   & \;[e_2,e^3]=e_2,       & \;[e_3,e^3]=-e^1,\cr
[e_1,e^4]=e_2,       & \;[e_2,e^4]=-(e_1+e^3),  & \;[e_3,e^4]=e^2. \cr
\end{array} \]
${\mathcal D}\ggo$ is $3$-step solvable, and observe that
\[\zg({\mathcal D}\ggo)=\text{span}\left\{e_4,e_3-e^1\right\},\qquad
[{\mathcal D}\ggo,{\mathcal
D}\ggo]=\text{span}\left\{e_2,e_3,e_4,e^1,e^2,e_1+e^3\right\},\] and
${\mathcal D}\ggo / \zg({\mathcal D}\ggo)$ is isomorphic to
$\RR^2\ltimes \RR^4$, where $
\RR^2=\text{span}\left\{x_1,x_2\right\}$ acts on $ \RR^4$ as
follows:
\[ \ad(x_1)=\begin{pmatrix}1&0& &  \\
  0& 1& &  \\
  & & -1&0  \\
  & &0 & -1 \end{pmatrix},\qquad
\ad(x_2)=\begin{pmatrix}0&-1& &  \\
  1& 0& &  \\
  & & 0&-1  \\
  & &1& 0 \end{pmatrix}.
\] The ad-invariant neutral metric $\langle \cdot\, ,\cdot \rangle$
on ${\mathcal D}\ggo$ is given by $\langle e_j,e^j \rangle=1$ for
$j=1,\ldots,4$, and the induced complex structure ${\mathcal J}$ on
${\mathcal D}\ggo$ is as follows:
\[ {\mathcal J}|_{\ggo}=J, \qquad \mathcal Je^1=e^2,\qquad \mathcal Je^3=-e^4.\]
\end{exam}

\medskip

\subsection{Compact semisimple bialgebras of complex type}
We recall Samelson's construction of a complex structure on a
compact semisimple even dimensional Lie algebra $\ggo$ \cite{S}.
Let $\hg$ be a maximal abelian subalgebra of $\ggo$. Then we have
the root space decomposition of $\ggo^{\CC}$ with respect to
$\hg^{\CC}$: \[ \ggo^{\CC}= \hg^{\CC} \oplus \sum _{\alpha \in
\Phi} \ggo ^{\CC}_{\alpha},\] where $\Phi$ is a finite subset of
$(\hg^{\CC})^*$ and \[\ggo^{\CC}_{\alpha}=\{\, x \in \ggo^{\CC}\,
:\, [h,x]=\alpha(h)x, \; \forall h \in  \hg^{\CC} \, \}\] are the
one dimensional root subspaces. Since $\hg$ is even dimensional,
one can choose a skew-symmetric endomorphism $I$ of $\hg$ with
respect to the Killing form such that $I^2=-$id. Samelson defines
a complex structure on $\ggo$ by considering a positive system
$\Phi^+$ of roots, which is a set $\Phi^+ \subset \Phi$ satisfying
\[ \Phi^+ \cap (-\Phi^+) = \emptyset, \qquad \Phi^+ \cup (-\Phi^+)
=\Phi, \qquad \alpha , \beta \in \Phi ^+ , \; \alpha +\beta \in
\Phi \Rightarrow \alpha +\beta \in \Phi ^+ .\] Setting \[ \mg =
\hg^{1,0} \oplus \sum _{\alpha \in \Phi^+}  \ggo ^{\CC}_{\alpha},
\] where $\hg^{1,0}$ is the eigenspace of $I$ of eigenvalue $i$,
it follows that $\mg$ is a (solvable) complex Lie subalgebra of
$\ggo^{\CC}$ which induces a complex structure $J$ on $\ggo$ such
that $\ggo^{1,0}=\mg$. This complex structure is skew-symmetric
with respect to the Killing form on $\ggo$, hence we obtain a Lie
bialgebra $(\ggo,\delta _J)$ of complex type. This is not
equivalent to the standard Lie bialgebra structure on $\ggo$
constructed in \cite{LW} (see also \cite{LR}). In fact, in the
standard case, the Lie algebra $\ggo^*$ is completely
solvable\footnote{Recall that a real solvable Lie algebra $\ug$ is
{\em completely solvable} when $\ad(x)$ has real eigenvalues for
all $x\in\ug$.}, whereas in the complex type case, it is not.

\smallskip

\begin{exam}
Let us consider $\ggo=\sug(2n+1)$ for $n\geq 2$, hence
$\ggo^\CC=\slg(2n+1,\CC)$. We can take
\[ \hg=\left\{\begin{pmatrix} it_1 & &\cr & \ddots & \cr & & it_{2n+1}\end{pmatrix} : t_j\in\RR,\,
t_1+\cdots+t_{2n+1}=0\right\}, \] and
\[ \hg^{1,0}=\left\{\begin{pmatrix} z_1 & &   \cr & \ddots &  \cr  & &  z_{2n+1}\end{pmatrix}
: \begin{array}{l} z_l\in\CC,   \;\; z_{1}-z_{2n+1}=iz_{n+1},   \\
  z_{r}=iz_{n+r},\,  r =2, \dots , n, \;\;
\sum_l z_l=0
\end{array} \right\}.
\] The corresponding complex structure $J$ on $\ggo$ as above is
determined by the complex subalgebra $\ggo^{1,0}$ of
$\slg(2n+1,\CC)$ given by $\ggo^{1,0}=\hg^{1,0}\oplus \ngo$ where
\[ \ngo=\{ \, \text{all } (2n+1)\times(2n+1)  \text{ strictly upper
triangular complex matrices} \, \}.\] According to
Theorem~\ref{bialg} and Proposition~\ref{prop1}, the Lie algebra
$\ggo^*$ associated to $(\ggo, \delta_J)$ is isomorphic to
$\ggo^{1,0}$.

On the other hand, we recall from \cite{LW} that if we consider the
standard Lie bialgebra structure on $\ggo$, then $\ggo^*$ is
isomorphic to the solvable Lie algebra
\[ \mathfrak{sb}(2n+1,\CC)= \left\{ \begin{array}{l} \text{all } (2n+1)\times(2n+1)  \text{ traceless upper
triangular} \\ \text{complex matrices with real diagonal entries}
\end{array}\right\}.
\] Note that $\mathfrak{sb}(2n+1,\CC)$ is not a complex subalgebra
of $\slg(2n+1,\CC)$. Therefore, our construction yields a Lie
bialgebra structure on $\sug(2n+1)$ which is not equivalent to the
standard one.
\end{exam}

\medskip

We exhibit next families of metric Lie algebras admitting
skew-symmetric complex structures where the ad-invariant metric has
signature $(2,n)$.

\subsection{Examples arising from Lorentzian ad-invariant metrics}\label{lor}
We recall from \cite{Me} the classification of Lie algebras
admitting ad-invariant metrics of signature $(1,n-1)$.

\begin{thm}\cite{Me}\label{med}
Each indecomposable non simple Lie algebra with an ad-invariant
Lorentzian metric is isomorphic to exactly one Lie algebra in the
family \[\osc({\underline{\lambda}})=\text{span}\{e_0, e_1,
\hdots, e_{2m+1}\}\] with $\underline{\lambda}=(\lambda_1,
\lambda_2, \hdots, \lambda_m)$, $1=\lambda_1\leq \lambda_2 \leq
\hdots \leq \lambda_m$, where \[ [e_{2i-1}, e_{2i}]=\lambda_i e_0,
\qquad i=1, \hdots, m,\] and the adjoint action of $e_{2m+1}$ on
$\text{span}\{e_1, e_2,\hdots , e_{2m}\}$ is given by
\begin{equation}\label{al}
 A_{\underline{\lambda}} : =
  \begin{pmatrix}
  0 & - \lambda_1 \\
  \lambda_1 & 0 \\
 & &  \ddots \\
  & & & 0 & -\lambda_m \\
  & & & \lambda_m & 0
  \end{pmatrix}.
  \end{equation}
  \end{thm}
The simply connected Lie group corresponding to
$\osc({\underline{\lambda}})$ is known as the {\em oscillator}
group. The $4$-dimensional Lie algebra $\osc(1)$ will be denoted
simply by $\osc$. Since the metric on any Lie algebra of the family
$\osc({\underline{\lambda}})$ is Lorentzian these Lie algebras do
not admit Lie bialgebra structures of complex type
(Lemma~\ref{sign}). However, trivial central extensions of them do
admit such structures.

\begin{thm}
The direct extensions $\osc({\underline{\lambda}})\times \RR^{1,1}$
are Lie bialgebras of complex type.
\end{thm}

\begin{proof}
Let $\osc({\underline{\lambda}})=\text{span}\{e_0, e_1, \hdots,
e_{2m+1}\}$ be  as above. An ad-invariant metric on
$\osc({\underline{\lambda}})$ is defined by the following non
trivial relations \[g(e_0, e_{2m+1})= g(e_{j}, e_{j})=1,\qquad
\text{ for  } j= 1, \hdots, 2m.\] If we consider $\RR^{1,1}=
\text{span} \{e_{2m+2}, e_{2m+3}\}$ equipped with the  metric
$g(e_{2m+2}, e_{2m+3})=1$ then we obtain on
$\osc({\underline{\lambda}})\times\RR^{1,1}$ an ad-invariant metric
of signature $(2,2m+2)$. The almost complex structure defined by
\[Je_0 = e_{2m+3}, \qquad Je_{2m+1}= e_{2m+2}, \qquad
Je_{2i-1}=e_{2i}, \qquad \text{ for  } i=1, \hdots, m,\] is
integrable and skew-symmetric with respect to $g$, proving the
assertion.
\end{proof}

\medskip

\subsection{Examples arising from ad-invariant metrics of signature $(2,n-2)$}\label{lor+1}
We recall from \cite{BK} the classification of the Lie algebras
with one dimensional centre admitting an ad-invariant metric of
signature $(2,n-2)$.

\begin{thm}\cite{BK}
Let $\ggo$ be an indecomposable Lie algebra with one dimensional
centre endowed with an ad-invariant metric of signature $(2,n-2)$.
Then $\ggo$ is isomorphic to one of the following Lie algebras:
$L_2(1,1), \; L_3(1,2), \; L_{2,\underline{\lambda}}(1, n-3)\,$ or
$\, L_{3, \underline{\lambda}}(1, n-3)$.
\end{thm}

We describe below the Lie algebras mentioned in the above theorem:
$L_2(1,1)$ was introduced in Example \ref{l2}, and for the
remaining cases we have

$\diamond\; L_3(1,2)=\text{span}\{e_0, e_1, e_2, e_3,e_4\}$ with
$[e_1,e_2]=e_0$ and the adjoint action of $e_4$ on the subspace
$\text{span}\{e_1,e_2, e_3\}$ is given by
\[L_3 = \begin{pmatrix}
 0 & 0 & 0\\
 1& 0 & 0 \\
 0 & 1 & 0
 \end{pmatrix}.\]

$\diamond\; L_{2,\underline{\lambda}}(1,n-3)=\text{span}\{e_0,
e_1, \hdots, e_{n-1}\}$ for $n=2m+2>5$ even and
$\underline{\lambda}=(\lambda_1, \ldots ,\lambda_m)$, $0<
\lambda_1\leq \cdots \leq \lambda_m,$ with $[e_1,e_2]=e_0$,
$[e_{2i-1}, e_{2i}]= \lambda_i e_0,\; 2\leq i \leq m$, and the
adjoint action of $e_{n-1}$ on $\text{span}\{e_1,\hdots,
e_{n-2}\}$ is given by \[
\begin{pmatrix}
 0 & 1 & \\
 1 & 0 & \\
 & & A_{\underline{\lambda}}
   \end{pmatrix},\]
with $A_{\underline{\lambda}}$ as in Theorem~\ref{med}.

$\diamond\; L_{3,\underline{\lambda}}(1,n-3)=\text{span}\{e_0, e_1,
\hdots, e_{n-1}\}$ for $n=2m+3>5$ odd, with $[e_1,e_2]=e_0$,
$[e_{2i}, e_{2i+1}]= \lambda_i e_0, \; 2\leq i \leq m$, and the
adjoint action of $e_{n-1}$ on $\text{span}\{e_1, \hdots, e_{n-2}\}$
is given by \[
\begin{pmatrix} L_3& 0 \\ 0 & A_{\underline{\lambda}}
\end{pmatrix},\]
with $\underline{\lambda}$ and $A_{\underline{\lambda}}$ as above.

\smallskip

\begin{thm}
Let $\ggo$ denote a Lie algebra as in the previous theorem. Then $
\ggo \times \RR^s$ is a Lie bialgebra of complex type, where $s=0$
(resp. $s=1$) if $\dim\ggo $ is even (resp. odd).
\end{thm}

\begin{proof}
We exhibit in each case a Hermitian structure $(J,g)$ where $g$  is
an ad-invariant metric. The proofs of the ad-invariance property and
the integrability of $J$ follow by standard computations.

\smallskip

$\diamond \;L_2(1,1):$ see Example~\ref{l2}.

\smallskip

$\diamond \; L_3(1,2)\times \RR e_5:$ \begin{gather*} g(e_0, e_4)
= -g(e_1, e_3) = g(e_2, e_2)=g(e_5,e_5)=1,\\
J e_0 = e_3, \qquad Je_1=e_4,\qquad J e_2=e_5. \end{gather*}

\smallskip

$\diamond \;L_{2,\underline{\lambda}}(1,n-3):$ \begin{gather*}
g(e_0, e_{n-1}) = -g(e_1, e_1) = g(e_i, e_i)= 1, \qquad   2\leq i\leq  n-2,\\
Je_0 = \frac1{\sqrt{2}}\, ({e_1 + e_2}), \qquad J e_{n-1} =
\frac1{\sqrt{2}}\, ({e_2 - e_1}),\qquad J e_{2i-1} = e_{2i},
\qquad 2\leq i\leq m.\end{gather*}

\smallskip

$\diamond \; L_{3,\underline{\lambda}}(1,n-3)\times \RR e_{n}:$
\begin{gather*} g(e_0, e_{n-1}) = -g(e_1, e_3) =  g(e_{2i},
e_{2i})=g(e_{2i+3}, e_{2i+3})=1, \qquad 1\leq i\leq  m,\\
J e_3 = e_0, \qquad Je_{n-1}=e_1,\qquad J e_{n}=e_2, \qquad J
e_{2i} = e_{2i+1},\qquad 2\leq i\leq m.\end{gather*}
\end{proof}

We can obtain Lie bialgebras of complex type in arbitrarily high
dimensions by applying an iterative doubling procedure. In fact,
it is a consequence of (1) in Corollary~\ref{bi-doble} that if
$(J,g)$ is a Hermitian structure on $\ggo$ such that $g$ is
ad-invariant, then ${\mathcal D}^k \ggo$ is a Lie bialgebra of
complex type of dimension $2^k m, \; m=\dim \ggo, \; k\geq 1$,
where $\mathcal D ^k \ggo$ is defined inductively by $\mathcal D
^1 \ggo =\mathcal D \ggo,\; \mathcal D^k \ggo =\mathcal
D\left(\mathcal D^{k-1} \ggo \right)$.

In a similar way, the second assertion in Corollary~\ref{bi-doble}
implies that if $J$ is a complex structure on any Lie algebra
$\ggo$, then $J$ induces a Hermitian structure on $T^{*}\ggo$ with
respect to the neutral ad-invariant metric, therefore, $\mathcal
D^k\left(T^{*} \ggo \right)$ is a Lie bialgebra of complex type
for $k\geq 1$.

We point out that the ad-invariant metric on either $\mathcal D ^k
\ggo$ or $\mathcal D^k\left(T^{*}\ggo \right)$ has neutral
signature for all $k\geq 1$. In the next section we will develop a
method to construct a Lie bialgebra of complex type with
ad-invariant metric of signature $(2p+2,2q+2)$ starting with one
whose corresponding metric has signature $(2p,2q)$.

\

\section{Construction of examples via double extensions}
The aim of this section is to obtain new examples of Lie
bialgebras of complex type starting with a Hermitian structure on
a lower dimensional Lie algebra by a double extension process. As
a consequence of our construction we show that all solvable metric
Lie algebras of dimension $4$ and $6$ with metrics of suitable
signature admit Lie bialgebra structures of complex type.

We begin by recalling the following  result (see, for instance,
\cite{MR}):

\begin{thm}
Let $(\dg,\langle \cdot \, ,\cdot\rangle)$ be an indecomposable
metric Lie algebra. Then one of the following conditions holds:
\begin{enumerate}
\item $\dg$ is simple;
\item $\dg$ is 1-dimensional;
\item $\dg$
is a double extension $\dg=\dg_\pi(\ggo,\hg)$, where $\hg$ is
1-dimensional or simple and $\pi:\hg \to \Der_{{\rm skew}}(\ggo
,\langle \cdot\, ,\cdot\rangle)$ is a representation of $\hg$ on
$\ggo$ by skew-symmetric derivations.
\end{enumerate}
\end{thm}

We will restrict ourselves to the solvable case, and therefore we
will consider only double extensions of a metric Lie algebra by
1-dimensional Lie algebras. Following the notation in \cite{BK}, we
give the description of such extensions. Let $\hg=\RR H$ and
$\hg^*=\RR \alpha$, with $\alpha(H)=1$. The double extensions of the
Lie algebra $(\ggo,\langle\cdot\, ,\cdot\rangle)$ by the one
dimensional Lie algebra $\hg$ are determined by skew-symmetric
derivations $A\in\Der_{\rm skew}(\ggo,\langle\cdot\,
,\cdot\rangle)$:
\[ \dg_A:=\dg_A(\ggo,\RR):=\RR \alpha\oplus\ggo\oplus\RR H\]
with Lie bracket
\[ \alpha\in\zg(\dg_A),\qquad   [x,y]_{\dg_A}=\langle Ax,y\rangle\alpha+[x,y]_{\ggo},\qquad [H,x]_{\dg_A}=Ax,\]
for all $x,y\in\ggo$. The ad-invariant metric on the double
extension $\dg_A$ is obtained extending the ad-invariant metric on
$\ggo$ via the single relation $\langle\alpha,H\rangle=1$.

We will show next that beginning with a metric Lie algebra
$(\ggo,\langle\cdot\, ,\cdot\rangle)$ equipped with a skew-symmetric
complex structure and certain skew-symmetric derivations, we can
produce a skew-symmetric complex structure on  double extensions of
$\ggo$.

\begin{thm}\label{ext}
Let $(J,\langle\cdot\, ,\cdot\rangle)$ be a Hermitian structure on
the  Lie algebra $\ggo$ with ad-invariant metric $\langle\cdot\,
,\cdot\rangle$, and let $\RR^{1,1}:=\text{span}\{x_0,y_0\}$ denote
the abelian pseudo-Euclidean Lie algebra with metric $\langle
x_0,y_0\rangle=1$. Denote by $\tilde{\ggo}$ the product metric Lie
algebra ${\tilde \ggo}=\RR^{1,1}\times\ggo$, and let ${\tilde A}$ be
a skew-symmetric derivation of ${\tilde\ggo}$. We denote by ${\JJ}$
the skew-symmetric almost complex structure on the double extension
$\dg_{\tilde A}=\dg_{\tilde A}({\tilde \ggo},\RR)$ defined by
\[ {\JJ}|_{\ggo}=J,\qquad {\JJ}x_0=\alpha,\qquad {\JJ}y_0=H,\qquad {\JJ}^2=-\text{id} . \]
Then ${\JJ}$ is integrable if and only if there exist $c\in\RR$,
${\overline y}\in\zg(\ggo)$, a 1-cocycle $k\in\ggo^*$ and a
skew-symmetric derivation $A$ of $(\ggo,\langle \cdot\,
,\cdot\rangle)$ such that
\[ {\tilde A}x_0=cx_0,\qquad {\tilde A}y_0=-cy_0+{\overline y},\qquad {\tilde A}x=k(x)x_0+Ax,\qquad AJ=JA, \]
where $k(x)=-\langle x,{\overline y}\rangle$ for all $x\in\ggo$.
\end{thm}

\begin{proof} In this proof we will denote $\dg:=\dg_{\tilde A}$, in order to simplify
notation. Let us first suppose that this almost complex structure
is integrable, and from this we will determine the form of the
derivation ${\tilde A}$. Recall that the integrability of ${\JJ}$
is equivalent to the vanishing of the Nijenhuis tensor $N_{\JJ}$,
where
\[ N_{\JJ}(x,y)={\JJ}[x,y]_{\dg}-[{\JJ}x,y]_{\dg}-[x,{\JJ}y]_{\dg}-
{\JJ}[{\JJ}x,{\JJ}y]_{\dg}, \qquad x,y\in\dg.\] From
$N_{\JJ}(x_0,y_0)=0$ we have that
\begin{align*}
0=N_{\JJ}(x_0,y_0) & ={\JJ}[x_0,y_0]_{\dg}-[{\JJ}x_0,y_0]_{\dg}-[x_0,{\JJ}y_0]_{\dg}-
{\JJ}[{\JJ}x_0,{\JJ}y_0]_{\dg}\\
& ={\JJ}\left(\langle {\tilde
A}x_0,y_0\rangle\alpha+[x_0,y_0]_{{\tilde
\ggo}}\right)-[\alpha,y_0]_{\dg}-[x_0,H]_{\dg}-
{\JJ}[\alpha,H]_{\dg}\\
& =-\langle {\tilde A}x_0,y_0\rangle x_0+{\tilde A}x_0.\\
\end{align*}
Hence, denoting $c:=\langle{\tilde A}x_0,y_0\rangle$, we obtain that
\[
{\tilde A}x_0=cx_0.\] Since ${\tilde A}$ is skew-symmetric, we
see that $\langle{\tilde A}y_0,y_0\rangle=0$, so that ${\tilde
A}y_0\in\left(\text{span}\{y_0\}\right)^{\perp}=\text{span}\{y_0\}\oplus\ggo$
and therefore ${\tilde A}y_0=\beta y_0+{\overline y}$ for some
$\beta\in\RR$ and ${\overline y}\in\ggo$. But, since ${\tilde A}$ is
skew-symmetric, we have
\[ c=\langle{\tilde A}x_0,y_0\rangle=-\langle x_0,{\tilde A}y_0\rangle=-\langle x_0,\beta y_0+{\overline y}\rangle=
-\beta\]
and using now that ${\tilde A}$ is a derivation of ${\tilde\ggo}$,
we obtain for $x\in\ggo$
\[ 0={\tilde A}[y_0,x]_{\tilde\ggo}=[{\tilde A}y_0,x]_{\tilde\ggo}+[y_0,{\tilde A}x_0]_{\tilde\ggo}=
[\beta y_0+{\overline y},x]_{\tilde\ggo}= [{\overline
y},x]_{\ggo},\] so that ${\overline y}\in\zg(\ggo)$. Thus, \[
{\tilde A}y_0=-cy_0+{\overline y},\qquad {\overline
y}\in\zg(\ggo).\]

Now take $x\in\ggo$ and observe that $\langle {\tilde
A}x,x_0\rangle=-\langle x,{\tilde A}x_0\rangle=-c\langle
x_0,x\rangle=0$, and hence ${\tilde
A}x\in\left(\text{span}\{x_0\}\right)^{\perp}=\text{span}\{x_0\}\oplus\ggo$,
so that there exist $k\in\ggo^*$ and $A\in\End(\ggo)$ such that \[
{\tilde A}x=k(x)x_0+Ax, \qquad\text{for all }x\in\ggo.\] Let us
determine $k$. For $x\in\ggo$, we compute $\langle {\tilde
A}y_0,x\rangle =\langle -cy_0+{\overline y},x\rangle=\langle
{\overline y},x\rangle$, but, on the other hand $\langle {\tilde
A}y_0,x\rangle =-\langle y_0,{\tilde A}x\rangle=-\langle
y_0,k(x)x_0+Ax\rangle=-k(x)$, so that \begin{equation}\label{k}
k(x)=-\langle x,{\overline y}\rangle,\qquad
x\in\ggo.\end{equation} In order to verify that $k$ is a
1-cocycle, we only have to show that $k([\ggo,\ggo]_{\ggo})=0$.
But this follows easily from \eqref{k} and the fact that
$\zg(\ggo)^{\perp}=[\ggo,\ggo]_{\ggo}$, since ${\overline
y}\in\zg(\ggo)$.

Now let us prove the statements regarding $A\in\End(\ggo)$. For
$x,y\in\ggo$, we have
\[\langle Ax,y\rangle=\langle {\tilde A}x-k(x)x_0,y\rangle=\langle{\tilde A}x,y\rangle=-\langle x,{\tilde A}y\rangle
=-\langle x,k(y)x_0+Ay\rangle=-\langle x,Ay\rangle,\] so that $A$ is
skew-symmetric. Next, we see that
\begin{align*}
A[x,y]_{\ggo} &={\tilde A}[x,y]_{\ggo}-k([x,y]_{\ggo})x_0=[H,[x,y]_{\ggo}]_{\dg}\\
              &=[H,[x,y]_{\dg}-\langle{\tilde A}x,y\rangle\alpha]_{\dg}=[H,[x,y]_{\dg}]_{\dg}\\
              &=-[x,[y,H]_{\dg}]_{\dg}-[y,[H,x]_{\dg}]_{\dg}\\
              &=[x,{\tilde A}y]_{\dg}-[y,{\tilde A}x]_{\dg}=[x,k(y)x_0+Ay]_{\dg}-[y,k(x)x_0+Ax]_{\dg}  \\
              &=-k(y)\underbrace{\langle {\tilde A}x_0,x\rangle}_{=0}\alpha+\langle{\tilde A}x,Ay\rangle\alpha+[x,Ay]_{\ggo}\\
              &\qquad +k(x)\underbrace{\langle{\tilde A}x_0,y\rangle}_{=0} \alpha-\langle{\tilde A}y,Ax\rangle\alpha-[y,Ax]_{\ggo}\\
              &=\langle k(x)x_0+Ax,Ay\rangle\alpha+[x,Ay]_{\ggo}-\langle k(y)x_0+Ay,Ax\rangle\alpha +[Ax,y]_{\ggo}\\
              &=[x,Ay]_{\ggo}+[Ax,y]_{\ggo},
\end{align*}
and hence $A$ is a derivation of $\ggo$. Finally, from $N_{\JJ}(y_0,x)=0$ we obtain that
\begin{align*}
0=N_{\JJ}(y_0,x) & ={\JJ}[y_0,x]_{\dg}-[H,x]_{\dg}-[y_0,Jx]_{\dg}-
{\JJ}[H,Jx]_{\dg}\\
& ={\JJ}\left(\langle {\tilde
A}y_0,x\rangle\alpha+[y_0,x]_{{\tilde\ggo}}\right)-{\tilde
A}x-\langle {\tilde
A}y_0,Jx\rangle\alpha-[y_0,Jx]_{\tilde\ggo}-{\JJ}{\tilde A}Jx\\
& =\langle y_0,{\tilde A}x\rangle x_0-(k(x)x_0+Ax)+\langle y_0,{\tilde A}Jx\rangle\alpha-{\JJ}\left(k(Jx)x_0+AJx\right)\\
& =k(x)x_0-k(x)x_0-Ax+k(Jx)\alpha-k(Jx)\alpha-JAJx\\
& =-Ax-JAJx
\end{align*}
for all $x\in\ggo$, so that $A=-JAJ$, or equivalently, $AJ=JA$.

Let us suppose now that ${\tilde A}$ is a skew-symmetric derivation
of $\left({\tilde\ggo},\langle\cdot \, ,\cdot \rangle\right)$
satisfying the conditions in the statement of the theorem. We have
to show that the Nijenhuis tensor $N_{\JJ}$ vanishes identically.
According to Lemma 2.2 in \cite{B}, we only have to verify that
$N_{\JJ}(x_0,y_0)=0,\;N_{\JJ}(x_0,x)=0,\; N_{\JJ}(y_0,x)=0$ and
$N_{\JJ}(x,y)=0$ for all $x,y\in\ggo$. From what has been done in
the first part of the proof, we see easily that
$N_{\JJ}(x_0,y_0)=0,\;N_{\JJ}(y_0,x)=0$ for all $x\in\ggo$; let us
consider the remaining cases. Take $x\in\ggo$ and compute
\begin{align*}
N_{\JJ}(x_0,x) & ={\JJ}[x_0,x]_{\dg}-[\alpha,x]_{\dg}-[x_0,{\JJ}x]_{\dg}-
{\JJ}[\alpha,{\JJ}x]_{\dg}\\
& ={\JJ}\left(\langle
{\tilde A}x_0,x\rangle\alpha+[x_0,x]_{{\tilde\ggo}}\right)-\left(\langle {\tilde A}x_0,Jx\rangle\alpha+
[x_0,Jx]_{\tilde\ggo}\right)\\
& =0,
\end{align*}
since ${\tilde A}x_0=cx_0$ is orthogonal to $\ggo$ and $x_0$ is
central in ${\tilde \ggo}$. Now, take $x,y\in\ggo$ and compute
\begin{align*}
N_{\JJ}(x,y) & ={\JJ}[x,y]_{\dg}-[Jx,y]_{\dg}-[x,Jy]_{\dg}-
{\JJ}[Jx,Jy]_{\dg}\\
& ={\JJ}\left(\langle {\tilde
A}x,y\rangle\alpha+[x,y]_{\ggo}\right)-\langle {\tilde
A}Jx,y\rangle\alpha-[Jx,y]_{\ggo}
-\langle{\tilde A}x,Jy\rangle\alpha-[x,Jy]_{\ggo}\\
& \qquad-{\JJ}\left(\langle {\tilde A}Jx,Jy\rangle\alpha+[Jx,Jy]_{\ggo}\right)\\
& =-\langle k(x)x_0+Ax,y\rangle x_0+J[x,y]_{\ggo}-\langle k(Jx)x_0+AJx,y\rangle\alpha-[Jx,y]_{\ggo}\\
& \qquad -\langle k(x)x_0+Ax,Jy\rangle\alpha -[x,Jy]_{\ggo}+\langle k(Jx)x_0+AJx,Jy\rangle x_0-J[Jx,Jy]_{\ggo}\\
& =-\langle Ax,y\rangle x_0-\langle JAx,y\rangle\alpha-\langle Ax,Jy\rangle\alpha+\langle JAx,Jy\rangle x_0+
N_J^{\ggo}(x,y)\\
& =-\langle Ax,y\rangle x_0+\langle Ax,Jy\rangle\alpha-\langle Ax,Jy\rangle\alpha+\langle Ax,y\rangle x_0+
N_J^{\ggo}(x,y)\\
& =0,
\end{align*}
since $J$ is a skew-symmetric integrable almost complex structure on
$\ggo$ which commutes with $A$. This  completes the proof.
\end{proof}

\smallskip

Let us compute next the associated Lie bracket $[\cdot,\cdot]_{\JJ}$
on $\dg_{\tilde A}$. Keeping the notation from the previous theorem,
it is easily verified that the following relations hold: for
$x,y\in\ggo$,
\[\begin{array}{ll}
[x,y]_{\JJ} = [x,y]_J, & [x_0,y_0]_{\JJ}=-cx_0,\cr [y_0,x]_{\JJ} =
k(x)x_0-k(Jx)\alpha+Ax, &  \cr [\alpha,y_0]_{\JJ}  = -c\alpha, &
[\alpha,H]_{\JJ}  =  cx_0, \cr [H,x]_{\JJ} =
k(x)\alpha+k(Jx)x_0+AJx, & [H,x_0]_{\JJ} = c\alpha.\\
\end{array}\]

From these equations we see that $\left(\dg_{\tilde A}\right)_{\JJ}=
\CC y_0 \ltimes  ( \CC x_0 \times \ggo _J)$ as complex Lie algebras,
where
\[ \ad_{\JJ}(y_0)= \begin{pmatrix} c& \vline &z_1& \cdots & z_n\\
\hline 0 & \vline & & & \\
 \vdots & \vline & & M+iN & \\
 0&\vline & &  & \end{pmatrix}
\]
with respect to an ordered $\CC$-basis $\{ x_0,v_1,\ldots , v_n\}$
of $\CC x_0 \times \ggo _J$. Here, $z_l=k(v_l)-ik(Jv_l)$, for
$l=1, \ldots, n$, and $A=\begin{pmatrix}M& -N\\N&M \end{pmatrix}$
in the ordered $\RR$-basis $\{v_1,\ldots , v_n, Jv_1,\ldots ,
Jv_n\}$ of $\ggo$.
\medskip

\subsection{Lie bialgebras of complex type in low dimensions} \label{low-dim}
In the next paragraphs we apply the method developed in
Theorem~\ref{ext}, starting with  $\ggo =\{ 0\}$ and $\ggo=\RR^2$,
to show that all solvable metric Lie algebras of dimension four
and six with metric of signature $(2r,2s)$ admit bialgebra
structures of complex type.

\medskip

$\ri$ We apply first Theorem~\ref{ext} in the case $\ggo=\{0\}$.
Hence ${\tilde\ggo}=\RR^{1,1}$ is spanned by $\{x_0,y_0\}$ with the
neutral metric $\langle x_0,y_0\rangle=1$ and the skew-symmetric
endomorphism ${\tilde A}$ is given  by ${\tilde A}x_0=cx_0,\,{\tilde
A}y_0=-cy_0$. In this case, the double extension $\dg_{\tilde
A}=\RR\alpha\oplus\RR^{1,1}\oplus\RR H$ has the following Lie
bracket relations
\[ [x_0,y_0]=c\alpha,\qquad [H,x_0]=cx_0,\qquad [H,y_0]=-cy_0, \]
the ad-invariant metric is given by $\langle
x_0,y_0\rangle=1=\langle\alpha,H\rangle$, and the skew-symmetric
complex structure ${\JJ}$ is determined by ${\JJ
}x_0=\alpha,\,{\JJ}y_0=H$.

If $c=0$, then $\dg_{\tilde A}$ is the pseudo-Euclidean abelian Lie
algebra with metric of signature $(2,2)$ and its canonical
skew-symmetric complex structure.

If $c\neq0$, then $(\dg_{\tilde A},{\JJ},\langle\cdot
\,,\cdot\rangle)$ is equivalent to $(L_2(1,1),J,g)$ from Example
~\ref{l2}, considering the following basis of $\dg_{\tilde A}$:
\[ e_1:=-c^{-1}H,\qquad e_2:=cy_0, \qquad e_3:=c^{-1}x_0, \qquad e_4:=-c\alpha.\]

$\rii$ Now we consider the case $\ggo=\RR^2$, the abelian
2-dimensional Euclidean Lie algebra. We take a basis $\{u,v\}$ of
$\ggo$ such that $\langle u,u\rangle=\langle v,v\rangle=1$ and
$Ju=v$, and hence ${\tilde\ggo}\cong \RR^{1,3}$ is spanned by
$\{x_0,y_0,u,v\}$ with the ad-invariant metric $\langle
x_0,y_0\rangle=\langle u,u\rangle=\langle v,v\rangle=1$. Any
skew-symmetric endomorphism ${\tilde A}$ of ${\tilde \ggo}$
satisfying the conditions of Theorem~\ref{ext} takes the following
form, in the ordered basis $\{x_0,y_0,u,v\}$:
\[ {\tilde A}=   \begin{pmatrix}
c&0&-a&-b\\0&-c&0&0\\0&a&0&-r\\0&b&r&0
\end{pmatrix}
\]
for some $a,b, c,r\in\RR$. For all values of these parameters, the
corresponding ad-invariant metric will have signature $(2,4)$. There
are two possible cases: the derivation ${\tilde A}$ is nilpotent or
not.

Suppose first that ${\tilde A}$ is nilpotent, that is, $c=r=0$.
This is equivalent to the Lie algebra $\dg_{\tilde
A}=\RR\alpha\oplus\RR^{1,3}\oplus\RR H$ being nilpotent (see
\cite{BK}). If $a=b=0$, we obtain that $\dg_{\tilde A}=\RR^{2,4}$.
On the other hand, if $a^2+b^2 \neq 0$,  we set
\begin{gather*} e_1:=(a^2+b^2)^{-1/4}H,\qquad e_2:=-(a^2+b^2)^{-1/4}y_0, \qquad e_3:=(a^2+b^2)^{-1/2}(-bu+av), \\
e_4:=(a^2+b^2)^{-1/2}(au+bv), \qquad e_5:=(a^2+b^2)^{1/4}x_0,
\qquad e_6:=(a^2+b^2)^{1/4}\alpha. \end{gather*} Thus for any
choice of such $a$ and $b$,  the Lie bracket on $\dg_{\tilde A}$
in terms of the basis $\{e_1,\ldots,e_6\}$ is given by
\begin{equation}\label{n1}
[e_1,e_2]=-e_4,\qquad [e_1,e_4]=-e_5,\qquad [e_2,e_4]=-e_6
,\end{equation} and the Hermitian structure $(J,g)$ on
$\dg_{\tilde A}$ becomes
\begin{gather*} Je_1=e_2,\qquad Je_3=-e_4, \qquad Je_5=e_6,\\
g(e_3,e_3)=g(e_4,e_4)=g(e_1,e_6)=-g(e_2,e_5)=1.\end{gather*} It can
be seen that $(\dg_{\tilde A},{\JJ},\langle\cdot \, ,\cdot \rangle)$
is equivalent to the metric Lie algebra with skew-symmetric complex
structure $(L_3(1,2)\times \RR,J,g)$ constructed in \S\ref{lor+1}.

\begin{rem}
The Lie algebra $L_3(1,2)\times \RR$ can be found in Salamon's
list of $6$-dimensional nilpotent Lie algebras admitting a complex
structure (see \cite{S}). In his notation, this Lie algebra
corresponds to $(0,0,0,12,14,24)$.
\end{rem}

\medskip

Assume next that ${\tilde A}$ is not nilpotent, so that
$c^2+r^2\neq 0$. Let us suppose that $a=b=0$, and hence
$k=0\in\ggo^*$ (we are using the notation from Theorem~\ref{ext}).
The Lie bracket on the double extension $\dg_{\tilde
A}=\RR\alpha\oplus\RR^{1,3}\oplus\RR H$ is given by
\begin{gather*} [x_0,y_0]=c\alpha, \qquad [H,x_0]=cx_0,\qquad [H,y_0]=-cy_0,\\
[H,u]=rv,\qquad [H,v]=-ru,\qquad [u,v]=r\alpha,\end{gather*} the
ad-invariant metric is obtained by adding the relation
$\langle\alpha,H\rangle=1$ to the ad-invariant metric on
${\tilde\ggo}$, and the skew-symmetric complex structure $\JJ$ is
defined as follows
\[ {\JJ}u=v,\qquad {\JJ}x_0=\alpha, \qquad {\JJ}y_0=H.\]
We have the following possibilities, depending on the conditions
satisfied by the parameters $c$ and $r$.

\medskip

$\diamond$ If $c=0,\,r\neq 0$, set
\[ e_1:=-r^{-1}H,\qquad e_2:=u,\qquad e_3:=v,\qquad e_4:=r\alpha,
\qquad e_5:=r^{-1}y_0, \qquad e_6:=rx_0.\] The Lie bracket on
$\dg_{\tilde A}$ in terms of the basis $\{e_1,\ldots,e_6\}$ is
given by
\[ [e_1,e_2]=-e_3,\qquad [e_1,e_3]=e_2, \qquad [e_2,e_3]=e_4,\]
and the Hermitian structure is
\begin{gather*}
Je_1=e_5,\qquad Je_2=e_3, \qquad Je_4=-e_6,\\
-g(e_1,e_4)=g(e_5,e_6)=g(e_2,e_2)=g(e_3,e_3)=1.
\end{gather*}
Therefore $(\dg_{\tilde A},{\JJ},\langle\cdot \,
,\cdot \rangle)$ is equivalent to the metric Lie algebra with
skew-symmetric complex structure
$(\mathfrak{osc}\times\RR^{1,1},J,g)$ constructed in \S\ref{lor}.

\smallskip

$\diamond$ If $c\neq 0,\,r=0$, set
\[ e_1:=-c^{-1}H,\qquad e_2:=c^{-1}y_0,\qquad e_3:=cx_0,\qquad e_4:=-c\alpha,
\qquad e_5:=u, \qquad e_6:=v.\] The Lie bracket on $\dg_{\tilde
A}$ in terms of the basis $\{e_1,\ldots,e_6\}$ is given by
\[ [e_1,e_2]=e_2,\qquad [e_1,e_3]=-e_3, \qquad [e_2,e_3]=e_4,\]
and the Hermitian structure is
\begin{gather*}
Je_1=e_2,\qquad Je_3=-e_4, \qquad Je_5=e_6,\\
g(e_1,e_4)=g(e_2,e_3)=g(e_5,e_5)=g(e_6,e_6)=1.
\end{gather*}
Hence, $(\dg_{\tilde A},{\JJ},\langle\cdot \,
,\cdot \rangle)$ is equivalent to the metric Lie algebra with
skew-symmetric complex structure $(L_2(1,1)\times\RR^2,J,g)$
constructed in \S\ref{lor+1}.
\smallskip

$\diamond$ If $c\neq 0,\,r\neq 0$, we will consider two cases,
according to the sign of $r/c$. If $r/c<0$, then we set
\[ e_1:=-c^{-1}H,\qquad e_2:=c^{-1}y_0,\qquad e_3:=cx_0,\qquad e_4:=u,
\qquad e_5:=v, \qquad e_6:=-c\alpha.\] On the other hand, if
$r/c>0$, we set
\[ e_1:=-c^{-1}H,\qquad e_2:=c^{-1}y_0,\qquad e_3:=cx_0,\qquad e_4:=v,
\qquad e_5:=u, \qquad e_6:=-c\alpha.\] In both cases,
$\{e_1,\ldots,e_6\}$ is a basis of $\dg_{\tilde A}$ such that its
Lie bracket is given by
\begin{equation} \label{again} \begin{array}{llll}
[e_1,e_2]=e_2, & [e_1,e_3]=-e_3, & [e_1,e_4]=\lambda e_5, &
[e_1,e_5]=-\lambda e_4, \cr [e_2,e_3]=e_6, & [e_4,e_5]=\lambda
e_6, & & \cr
\end{array}\end{equation}
for $\lambda=|r/c|$ and it can be shown that  $\dg_{\tilde A}$ is
isomorphic to $\ggo_{\lambda}:=L_{2,\lambda}(1,3)$ from
\S\ref{lor+1}. The ad-invariant metric $g$ is the same for both
cases and it is given by
\[ g(e_1,e_6)=g(e_2,e_3)=g(e_4,e_4)=g(e_5,e_5)=1, \]
whereas the skew-symmetric complex structures $J_1$ and $J_2$
corresponding to each case are given by
\[ J_1e_1=e_2,\qquad J_1e_3=-e_6, \qquad J_1e_4=e_5, \]
for $r/c<0$ and
\[ J_2e_1=e_2,\qquad J_2e_3=-e_6, \qquad J_2e_4=-e_5, \]
for $r/c>0$.  Recall that for any $\mu\in\CC$, the $3$-dimensional
complex Lie algebra $\rg_{3,\,\mu}(\CC)$ has a basis
$\{x_0,x_1,x_2\}$ with Lie brackets given by:
$[x_0,x_1]=x_1,\;[x_0,x_2]=\mu x_2$, and moreover,
$\rg_{3,\,\mu}(\CC)$ is isomorphic to $\rg_{3,\,1/{\mu}}(\CC)$ for
$\mu \neq 0$ (see \cite{GOV}). It turns out that
$\left(\ggo_{\lambda}\right)_{J_1}\cong \rg_{3,\,-i\lambda}(\CC)$
and $\left(\ggo_{\lambda}\right)_{J_2}\cong \rg_{3,\,i\lambda}(\CC)$
as complex Lie algebras.

\begin{rems}\

\smallskip

\noindent $\ri$ Since $\left(\ggo_{\lambda}\right)_{J_1}$ and
$\left(\ggo_{\lambda}\right)_{J_2}$ are not isomorphic for
$\lambda\neq 1$, $J_1$ and $J_2$ are not equivalent and the
corresponding bialgebras $(\ggo,\delta_{J_i}), \; i=1,2,$ are not
isomorphic. When $\lambda =1$, it can be shown that $J_1$ is still
not equivalent to $J_2$ and the corresponding bialgebras are not
isomorphic, even though, in this case,
$\left(\ggo_{1}\right)_{J_1}\cong\left(\ggo_{1}\right)_{J_2}$.
Therefore, there exist two non isomorphic Lie bialgebra structures
on $\ggo_1$ such that they induce isomorphic Lie algebra
structures on $\ggo_1^*$ (compare with Theorem 1.12 in \cite{LW}).

\smallskip

\noindent $\rii$ Concerning the behaviour of the Lie bialgebras
$(\ggo _{\lambda},\delta_{J_i} )$  for $i=1, 2$ as $\lambda$
approaches $0$ or $\infty$, we point out that when $\lambda \to
0$, we obtain in both cases ($i=1,2$) the Lie bialgebra
$L_2(1,1)\times\RR^2$  with the structure constructed  above. On
the other hand, if $\lambda \to \infty$, then  the resulting Lie
bialgebra is $\osc \times \RR^{1,1}$. Therefore,
$(\ggo_{\lambda},\delta_{J_1} )$ and $(\ggo_{\lambda},\delta_{J_2}
)$ are non-isomorphic Lie bialgebras which converge to the same
bialgebra when $\lambda\to 0$ or $\lambda \to \infty$.
\end{rems}

\medskip

In \cite{BK}, the indecomposable metric Lie algebras of dimension
$\leq 6$ have been classified. From this classification we can
determine all (not necessarily indecomposable) solvable metric Lie
algebras with metric of signature $(2r,2s)$ in these dimensions.
We list below the non-abelian ones:
\begin{equation}\label{list}
L_{2}(1,1),\quad  L_3(1,2)\times\RR, \quad L_{2}(1,1)\times \RR^2, \quad \osc
\times\RR^{1,1},\quad L_{2,\lambda}(1,3) \; (\lambda
>0).\end{equation}
Therefore, the results in this section can be summarized as
follows:

\begin{thm}\label{teo-low}
Every solvable metric Lie algebra $(\ggo, g),\;\dim\ggo \leq 6$,
such that $g$ has signature $(2r,2s)$ admits a skew-symmetric
complex structure and hence, a Lie bialgebra structure of complex
type.
\end{thm}

\smallskip

According to Theorem~\ref{bialg}, for any Lie bialgebra of complex
type $(\ggo,\delta_J)$, its dual $\ggo^*$ inherits a complex Lie
algebra structure. Table~\ref{comm} exhibits $\ggo^*$ when $\ggo$ is
one of the algebras in \eqref{list} with the Hermitian structures
considered above (compare with the Remark~(ii) after
Proposition~\ref{prop1}). The complex Lie algebras $\aff(\CC)$ and
$\rg_{3,\pm i\lambda}(\CC)$ in Table~\ref{comm} have already been
introduced, and $\hg_3(\CC)$ denotes the $3$-dimensional complex
Heisenberg Lie algebra.

\begin{table}
\begin{center}
\begin{tabular}{|l|l|}\hline
$\ggo$ & $\ggo^* $    \\
\hline 
$ L_{2}(1,1)$ & $\aff(\CC)$\\
$ L_3(1,2)\times\RR$ &
$\hg_3(\CC)$ \\
$L_{2}(1,1)\times \RR^2$ &
$\aff(\CC)\times\CC$\\
$\osc\times\RR^{1,1}$ & $\aff(\CC)\times\CC$ \\
$L_{2,\lambda}(1,3) $ & $\begin{cases} \rg_{3,-i\lambda}(\CC) \text{ for the complex structure } J_1  \\ 
\rg_{3,i\lambda}(\CC) \text{ for the complex structure } J_2 \\
\end{cases} $\\
\hline
\end{tabular}
\bigskip
\caption{}\label{comm}
\end{center}
\end{table}

\medskip

\begin{rem} We point out that Theorem~\ref{teo-low} is no longer true
in dimension $8$. In fact, it follows from \cite{ABDF} that the
cotangent algebra $(T^* \mathfrak r_4 ,\langle \, \cdot , \cdot
\rangle)$ endowed with its natural neutral ad-invariant  metric
$\langle \, \cdot , \cdot \rangle$  of signature $(4,4)$ does not
admit any skew-symmetric complex structure, where the Lie algebra
$\mathfrak r_4 =\text{span}\{e_0,e_1,e_2, e_3\}$ is defined by \[
[e_0,e_1]=e_1, \qquad [e_0, e_2]=e_1+e_2, \qquad
[e_0,e_3]=e_2+e_3.\]
\end{rem}

\

\section{Symplectic foliations on some Poisson Lie groups of complex type}\label{examples}

In this section we exhibit some examples of low dimensional simply
connected Poisson Lie groups of complex type. The underlying Lie
groups are either solvable or nilpotent and thus they are
diffeomorphic to Euclidean space. We also determine the
corresponding symplectic foliations in some cases.

In what follows we will denote
$\partial_i:=\frac{\partial}{\partial x_i}$ and $\partial
_{i,j}:=\frac{\partial}{\partial x_i}\wedge
\frac{\partial}{\partial x_j}$. Also, for any function $h\in
C^\infty(\RR^n)$ we set $h_i:=\frac{\partial h}{\partial x_i}$ and
the Hamiltonian vector field $X_h$ associated to $h$ is defined by
$X_h(g)=\{h,g\}$. If $\Pi$ is the Poisson structure corresponding
to $\{\cdot , \cdot \}$, the characteristic distribution
associated to $\Pi$ is given by
\[ {\mathcal C}_x=\{\left. X_h\right|_x:h\in C^\infty(\RR^n)\}\subset T_x\RR^n.\]

\begin{exam}
We consider the Lie algebra $\ggo:=L_2(1,1)$ from
Example~\ref{l2}. Let $G$ denote the simply connected Lie group
with Lie algebra $\ggo$; it is diffeomorphic to $\RR^4$ and it can
be realized as the following matrix group:
\[ G=\left\{\begin{pmatrix} 1 & x_2 & x_4\cr 0 & {\rm e}^{-x_1} & x_3\cr 0 & 0 & 1\cr\end{pmatrix}\,:
\,x_i\in\RR\right\}.\] Identifying the matrix  above with the point
$(x_1,x_2,x_3,x_4)\in\RR^4$, then $\RR^4$ acquires a group structure
given by
\[ (x_1,x_2,x_3,x_4)*(y_1,y_2,y_3,y_4)=(x_1+y_1,x_2\e^{-y_1}+y_2,y_3\e^{-x_1}+x_3,x_4+y_4+x_2y_3).\]
The r-matrix $J:\ggo\rightarrow\ggo$ can be considered as an
element $R\in\alt^2\ggo$ via the ad-invariant metric
$g:\ggo\rightarrow\ggo^*$, by setting $R:=J\circ
g^{-1}:\ggo^*\rightarrow\ggo$. In this case, one can easily see
that $R=e_1\wedge e_3-e_2\wedge e_4$, where $e_1,\ldots,e_4$ are
as in \eqref{cor-l2}. According to the remark following
Proposition~\ref{r-bialg}, the corresponding multiplicative
Poisson tensor $\Pi_R$ determined by $R$ making $G$ a Poisson Lie
group is given by
$\Pi_R=\,\stackrel{\leftarrow}{R}-\stackrel{\rightarrow}{R}$.
Performing standard computations, we arrive at the following
expression for $\Pi_R$:
\[ \Pi_R=(\e^{-x_1}-1){\partial}_{1,3}+x_2{\partial}_{1,4}-x_2\e^{-x_1}{\partial}_{2,3}+
(\e^{-x_1}-x_2^2-1){\partial}_{2,4}. \] Equivalently, the Poisson
bracket of two functions $f,g\in C^{\infty}(\RR^4)$ is given by
\begin{multline*} \{f,g\}=(\e^{-x_1}-1)(f_1g_3-f_3g_1)+x_2(f_1g_4-f_4g_1)\\-x_2\e^{-x_1}(f_2g_3-f_3g_2)
+(\e^{-x_1}-x_2^2-1)(f_2g_4-f_4g_2).\end{multline*}
From this, we obtain that the Hamiltonian vector field $X_f$
associated to $f\in C^\infty(\RR^4)$ can be expressed as
\begin{multline*}
X_f=-\left((\e^{-x_1}-1)f_3+x_2f_4\right){\partial}_1+\left(x_2\e^{-x_1}f_3-(\e^{-x_1}-x_2^2-1)f_4\right){\partial}_2
\\
+\left((\e^{-x_1}-1)f_1-x_2\e^{-x_1}f_2\right){\partial}_3+
\left(x_2f_1+(\e^{-x_1}-x_2^2-1)f_2\right){\partial}_4.
\end{multline*}

Let us consider now any point $p=(x_1,x_2,x_3,x_4)\in\RR^4$. It is
easy to see that $\Pi_R(p)=0$ if and only if $x_1=x_2=0$, whereas if
$x_1^2+x_2^2\neq 0$, then $\Pi_R(p)$ has rank $4$, where we consider
$\Pi_R(p)$ as a linear transformation $\Pi_R(p):T^*_p\RR^4\lra
T_p\RR^4$. Consequently, the symplectic foliation on $\RR^4$
determined by $\Pi_R$ can be described in the following way: each
point $\{(0,0,x_3,x_4)\}$ is a $0$-dimensional leaf, there are no
$2$-dimensional leaves, and the open set
\[ \mathcal U=\{(x_1,x_2,x_3,x_4)\in\RR^4:x_1^2+x_2^2\neq 0\} \]
is the remaining $4$-dimensional symplectic leaf. The corresponding
symplectic form on $\mathcal U$ $\omega=(\Pi_R|_{\mathcal U})^{-1}$
is given by
\begin{multline*} \omega=\frac{1}{\Delta}\left\{-(\e^{-x_1}-x_2^2-1)\dif x_1\wedge\dif
x_3-x_2\e^{-x_1}\dif x_1\wedge\dif x_4 \right. \\ \left. +x_2\dif x_2\wedge\dif
x_3-(\e^{-x_1}-1)\dif x_2\wedge\dif x_4\right\},\end{multline*} where
$\Delta=(\e^{-x_1}-1)^2+x_2^2$.

Let \[ \Pi_R^{(1)}=-x_1
{\partial}_{1,3}+x_2{\partial}_{1,4}-x_2{\partial}_{2,3}-
x_1{\partial}_{2,4}
\] be the linear part of $\Pi_R$ in a neighbourhood of $(0,0,x_3,x_4)$.
It follows that the Lie algebra determined by $\Pi_R^{(1)}$ is
isomorphic to $\aff(\CC)$. Therefore, according to \cite{DZ},
$\Pi_R$ is analytically linearizable.
\end{exam}

\medskip

\begin{exam}
Let $\ngo:=L_3(1,2)\times\RR$ be the nilpotent Lie algebra
considered in \S\ref{low-dim}, and denote by ${\mathcal N}$ the
corresponding simply connected Lie group, which is diffeomorphic
to $\RR^6$ and can be realized as the following matrix group:
\[ {\mathcal N}=\left\{\begin{pmatrix} 1 & 0 & -\frac{2}{3}x_2^2 & 0 & \frac{1}{6}x_1x_2+\frac{1}{3}x_4 &
-\frac{1}{6}x_2^2 & x_6 \cr
0 & 1 & -\frac{2}{3}x_1 & 0 & \frac{1}{6}x_1^2 &
-\frac{1}{6}x_1x_2+\frac{1}{3}x_4 & x_5 \cr 0 & 0 & 1 & 0 &
-\frac{1}{2}x_1 & \frac{1}{2}x_2 & x_4 \cr 0 & 0 & 0 & 1 & 0 & 0 &
x_3 \cr 0 & 0 & 0 & 0 & 1 & 0 & x_2 \cr 0 & 0 & 0 & 0 & 0 & 1 &
x_1 \cr 0 & 0 & 0 & 0 & 0 & 0 & 1\cr
\end{pmatrix}\,:\,x_i\in\RR\right\}.\]
Identifying the matrix  above with the point
$(x_1,\ldots,x_6)\in\RR^6$, then $\RR^6$ acquires a group structure
given by
\[ \begin{pmatrix} x_1\cr x_2\cr x_3\cr x_4\cr x_5\cr x_6\end{pmatrix}*
\begin{pmatrix} y_1\cr y_2\cr y_3\cr y_4\cr y_5\cr y_6\end{pmatrix}=
\begin{pmatrix} x_1+y_1 \cr x_2+y_2 \cr x_3+y_3 \cr x_4+y_4-\frac{1}{2}(x_1y_2-x_2y_1)\cr
x_5+y_5+\frac{1}{6}x_1^2y_2-\frac{2}{3}x_1y_4+(-\frac{1}{6}x_1x_2+\frac{1}{3}x_4)y_1\cr
x_6+y_6-\frac{1}{6}x_2^2y_1-\frac{2}{3}x_2y_4+(\frac{1}{6}x_1x_2+\frac{1}{3}x_4)y_2\end{pmatrix}.\]
The r-matrix $J:\ngo\rightarrow\ngo$ determines a bivector
$R\in\alt^2\ngo$ by setting $R:=J\circ
g^{-1}:\ngo^*\rightarrow\ngo$. In this case, one can easily see
that $R=e_4\wedge e_3+e_5\wedge e_1+e_6\wedge e_2$, where $e_1,
\dots ,e_6$ are as in \eqref{n1}, and the multiplicative Poisson
tensor
$\Pi_R=\,\stackrel{\leftarrow}{R}-\stackrel{\rightarrow}{R}$
determined by $R$ making ${\mathcal N}$ a Poisson Lie group is
given by:
\[ \Pi_R=x_1{\partial}_{3,5}+x_2{\partial}_{3,6}-x_2{\partial}_{4,5}
+x_1{\partial}_{4,6}-\frac{1}{6}(x_1^2+x_2^2){\partial}_{5,6}. \]
Note that $\Pi_R(x_1,\ldots,x_6)=0$ if and only if $x_1=x_2=0$. The
Poisson bracket of two functions $f,g\in C^{\infty}(\RR^6)$ induced
by $\Pi_R$ is given by
\begin{multline*}  \{f,g\}=x_1(f_3g_5-f_5g_3)+x_2(f_3g_6-f_6g_6)-x_2(f_4g_5-f_5g_4)\\+x_1(f_4g_6-f_6g_4)
-\frac{1}{6}(x_1^2+x_2^2)(f_5g_6-f_6g_5).\end{multline*} Therefore,
the Hamiltonian vector field $X_f$ associated to $f\in
C^\infty(\RR^6)$ can be expressed as
\begin{multline*} X_f=\left(-x_1f_5-x_2f_6\right){\partial}_3+\left(x_2f_5-x_1f_6\right){\partial}_4\\
+\left(x_1f_3-x_2f_4+\frac{1}{6}(x_1^2+x_2^2)f_6\right){\partial}_5
+\left(x_2f_3+x_1f_4-\frac{1}{6}(x_1^2+x_2^2)f_5\right){\partial}_6.
\end{multline*}
The symplectic foliation on $\RR^6$ corresponding to $\Pi_R$ is the
singular foliation associated to the characteristic distribution
${\mathcal C}_x$, and in what follows we are to determine its
symplectic leaves. Clearly, ${\mathcal C}_x=\{0\}$ if and only if
$x_1=x_2=0$. Let us denote by $p_j,\,j=1,\ldots,6$, the projections
$p_j(x_1,\ldots,x_6)=x_j$. Then we have $X_{p_1}=X_{p_2}=0$ and also
\begin{align*}
X_{p_3} & =x_1{\partial}_5+x_2{\partial}_6,\\
X_{p_4} & =-x_2{\partial}_5+x_1{\partial}_6,\\
X_{p_5} &
=-x_1{\partial}_3+x_2{\partial}_4-\frac{1}{6}\left(x_1^2+x_2^2\right)
{\partial}_6,\\
X_{p_6} &
=-x_2{\partial}_3-x_1{\partial}_4+\frac{1}{6}\left(x_1^2+x_2^2\right)
{\partial}_5.
\end{align*}
Observe now that if $x_1^2+x_2^2\neq 0$, then we have that
$\{\left. X_{p_3}\right|_x,\left. X_{p_4}\right|_x,\left.
X_{p_5}\right|_x,\left. X_{p_6}\right|_x\}$ is a basis of
${\mathcal C}_x$. Moreover, from the expressions obtained for
$X_{p_j}(x)$ above, we can easily see that
$\left\{\left.{\partial}_3\right|_x,\left.{\partial}_4\right|_x,\left.{\partial}_5\right|_x,
\left.{\partial}_6\right|_x \right\}$ is another basis for
${\mathcal C}_x$. From all these considerations we obtain the
following description of the symplectic leaves:

$\bullet$ the $0$-dimensional leaves are the points
$\{(0,0,c_3,c_4,c_5,c_6)\}$;

$\bullet$ there are no $2$-dimensional
leaves;

$\bullet$ the $4$-dimensional leaves are the affine subspaces
\[{\mathcal F}_{c_1,c_2}=\{(c_1,c_2,x_3,x_4,x_5,x_6):
x_3,\ldots,x_6\in\RR\},\] where $c_1,c_2$ are real constants such
that $c_1^2+c_2^2\neq 0$. The symplectic form $\omega_{c_1,c_2}$
on ${\mathcal F}_{c_1,c_2}$ induced by $\Pi_C$ is given by
\begin{multline*}
\omega_{c_1,c_2}  =  -\frac{1}{6}\dif x_3\wedge\dif x_4 -(c_1^2+c_2^2)^{-1}(c_1\dif x_3\wedge\dif x_5 \\
+c_2\dif x_3\wedge\dif x_6-c_2\dif x_4\wedge\dif x_5+c_1\dif
x_4\wedge\dif x_6), \end{multline*} where $(x_3,\ldots,x_6)$ are
global coordinates on ${\mathcal F}_{c_1,c_2}$.

\smallskip

It follows from Proposition~3.3 in \cite{C-A} that $\Pi_R$ is
analytically linearizable. Let \[
\Pi_R^{(1)}=x_1{\partial}_{3,5}+x_2{\partial}_{3,6}-x_2{\partial}_{4,5}
+x_1{\partial}_{4,6}
\] be the linear part of $\Pi_R$ in a neighbourhood of $(0,0,x_3,x_4,x_5,x_6)$.
Observe that the Lie algebra determined by $\Pi_R^{(1)}$ is
isomorphic to the complex Heisenberg Lie algebra $\hg _3(\CC)$,
considered as a real Lie algebra.
\end{exam}

\begin{exam} We consider the Lie group $G_{\lambda}$ with Lie
algebra $\ggo_{\lambda}:=L_{2, \lambda}(1,3)$ from \S\ref{low-dim},
for $\lambda>0$. This group has the following matrix realization:
\[
G_{\lambda}=\left\{  \left(
                         \begin{array}{cccccc}
                           1 & x_2 & 0 & x_4 \sqrt{\frac{\lambda}2} & x_5 \sqrt{\frac{\lambda}2} & x_6 \\
                           0 & \e^{-x_1} & 0 & 0 & 0 & x_3 \\
                           0 & 0 & 1 & 0 & 0 & 0 \\
                           0 & 0 & 0 & \cos (\lambda x_1) & -\sin(\lambda x_1) &
                           \sqrt{\frac{\lambda}2} \,\left(x_4 \sin(\lambda x_1)+x_5\cos (\lambda x_1)\right) \\
                           0 & 0 & 0 & \sin(\lambda x_1) & \cos (\lambda x_1) &
                           \sqrt{\frac{\lambda}2} \, \left(-x_4 \cos(\lambda x_1)+x_5\sin (\lambda x_1)\right)\\
                           0 & 0 & 0 & 0 & 0 & 1 \\
                         \end{array}
                       \right)
    \right\},
\]
for $x_1, \dots , x_6 \in \RR$. This induces the following
multiplication on $\RR^6$:
\[ \begin{pmatrix} x_1\cr x_2\cr x_3\cr x_4\cr x_5\cr x_6\end{pmatrix}*
\begin{pmatrix} y_1\cr y_2\cr y_3\cr y_4\cr y_5\cr y_6\end{pmatrix}=
\begin{pmatrix} x_1+y_1 \cr y_2+x_2 \e^{-y_1} \cr x_3+y_3 \e^{-x_1}
\cr y_4+x_4 \cos(\lambda y_1)+x_5\sin(\lambda y_1)\cr y_5 - x_4
\sin(\lambda y_1)+x_5\cos(\lambda y_1) \cr x_6+y_6 +x_2 y_3+
\frac{\lambda}2 \sin(\lambda y_1)(x_4y_4+x_5y_5) +\frac{\lambda}2
\cos(\lambda y_1)(x_4y_5-x_5y_4)
\end{pmatrix}.
\]
Recall from \S\ref{low-dim} the two inequivalent Hermitian
structures $(J_1,g)$ and $(J_2,g)$ on $\ggo_{\lambda}$. The
r-matrices $J_i: \ggo_{\lambda}\rightarrow \ggo_{\lambda}$
determine a bivector $R_i\in\alt^2 \ggo_{\lambda}$ by setting
$R_i:=J_i\circ g^{-1}:\ggo_{\lambda}^*\rightarrow \ggo_{\lambda}$,
for $i=1,2$. One can easily see that
\begin{gather*}
R_1=e_1\wedge e_3-e_2\wedge e_6+e_4\wedge e_5,\\
R_2=e_1\wedge e_3-e_2\wedge e_6-e_4\wedge e_5,
\end{gather*}
where $e_1,\ldots,e_6$ are as in \eqref{again}, and the multiplicative Poisson tensors
$\Pi_i^\lambda=\,\stackrel{\leftarrow}{R_i}-\stackrel{\rightarrow}{R_i}$
determined by $R_i$ are given by:
\begin{multline*} \Pi_1^\lambda=(\e^{-x_1}-1)\partial _{1,3}+ x_2
 \partial _{1,6}-x_2 \e ^{-x_1} \partial _{2,3}
+(\e^{-x_1}-x_2^2-1)  \partial _{2,6} \\ -\lambda x_5 \e^{-x_1}
 \partial _{3,4}  +\lambda x_4 \e^{-x_1} \partial _{3,5}   +\lambda (x_4 +x_2
x_5) \partial _{4,6}  +\lambda (x_5-x_2 x_4)
\partial _{5,6},
\end{multline*}
\begin{multline*} \Pi_2^\lambda=(\e^{-x_1}-1)\partial _{1,3}+ x_2
 \partial _{1,6}-x_2 \e ^{-x_1} \partial _{2,3}
+(\e^{-x_1}-x_2^2-1)  \partial _{2,6} \\ -\lambda x_5 \e^{-x_1}
 \partial _{3,4}  +\lambda x_4 \e^{-x_1} \partial _{3,5}  -\lambda (x_4 -x_2 x_5)
 \partial_{4,6}  -\lambda (x_5+x_2 x_4)\partial _{5,6}.
\end{multline*}

Note that $\Pi_i^\lambda(x_1,\ldots,x_6)=0$ if and only if
$x_1=x_2=x_4=x_5=0$. The characteristic distribution ${\mathcal
C}_i^\lambda$ associated to $\Pi_i^\lambda$ satisfies $({\mathcal C}_i^\lambda)_x=\{0\}$
if and only if $x_1=x_2=x_4=x_5=0$. For the remaining points,
$\dim ({\mathcal C}_i^\lambda)_x=4$ and \[({\mathcal C}_i^\lambda)_x=
\text{span}\left\{ \left. V\right|_x,\left. W_i\right|_x,
\left.{\partial}_3\right|_x , \left.{\partial}_6\right|_x
\right\},\] where $V$ and $W_i$ are the following vector fields on
$\RR^6$:
\begin{gather*} V= -\left(\e^{-x_1}-1\right){\partial}_1 +
x_2\e^{-x_1}{\partial}_2 -\lambda x_5 \e^{-x_1}{\partial}_4 +\lambda
x_4 \e^{-x_1}{\partial}_5,
\\ W_1=-x_2 {\partial}_1
-\left(\e^{-x_1}-x_2^2-1\right){\partial}_2
-\lambda(x_4+x_2x_5){\partial}_4 -\lambda (x_5-x_2x_4){\partial}_5,\\
W_2=-x_2 {\partial}_1-\left(\e^{-x_1}-x_2^2-1\right){\partial}_2
+\lambda(x_4-x_2x_5){\partial}_4 +\lambda (x_5+x_2x_4){\partial}_5.
\end{gather*}
Note that $[V,W_i]=0$ for $i=1,2$. We conclude that in both
cases the $0$-dimensional leaves are the points
$\{(0,0,c_3,0,0,c_6)\}$ and the open set $\mathcal U=\{ x \in \RR^6
: x_1^2+x_2^2+x_4^2+x_5^2\neq 0 \} $ is foliated by $4$-dimensional
symplectic manifolds, the integral manifolds of the above
distributions. Consider the disjoint union $\mathcal U=\mathcal
U_1\cup \mathcal U_2 \cup \mathcal U_3$, where $\mathcal U_i,\,
i=1,2,3$, are the sets
\begin{eqnarray*}
\mathcal U_1 & = & \{x \in \RR^6 : x_1^2+x_2^2\neq 0,\, x_4=x_5=0 \},\\
\mathcal U_2 & = & \{x \in \RR^6 : x_4^2+x_5^2\neq 0,\, x_1=x_2=0 \},\\
\mathcal U_3 & = & \{x \in \RR^6 : x_1^2+x_2^2\neq 0,\,
x_4^2+x_5^2\neq 0 \}.
\end{eqnarray*}
The set $\mathcal U_i$, $i=1$ or $2$, has four connected
components diffeomorphic to $\RR^4$, and each of them is a leaf of
both distributions, since $({\mathcal C}_1^\lambda)_x=({\mathcal
C}_2^\lambda)_x$ for $x\in \mathcal U_1\cup \mathcal U_2$. On the
other hand, these distributions satisfy $({\mathcal
C}_1^\lambda)_x \neq ({\mathcal C}_2^\lambda)_x$ for each $x\in
\mathcal U_3$, and therefore, their integral manifolds through
these points do not coincide. Nevertheless, the integral manifolds
of ${\mathcal C}_1^\lambda$ and ${\mathcal C}_2^\lambda$ through
$c=(c_1,\ldots,c_6)\in \mathcal U_3$ are diffeomorphic to $S_c
\times \RR^2$, where $S_c$ is a surface contained in $\{x \in
\mathcal U_3 : x_3=c_3,\,x_6=c_6\}$.

Let \begin{gather*}
\left(\Pi_1^\lambda\right)^{(1)}=-x_1{\partial}_{1,3}+x_2{\partial}_{1,6}-x_2{\partial}_{
2,3} -x_1{\partial}_{2,6}
 -\lambda x_5{\partial}_{3,4}+\lambda x_4{\partial}_{3,5}+\lambda
x_4{\partial}_{4,6} +\lambda x_5{\partial}_{5,6} \\
\left(\Pi_2^\lambda\right)^{(1)}=-x_1{\partial}_{1,3}+x_2{\partial}_{1,6}-x_2{\partial}_{
2,3} -x_1{\partial}_{2,6}
 -\lambda x_5{\partial}_{3,4}+\lambda x_4{\partial}_{3,5}-\lambda
x_4{\partial}_{4,6} -\lambda x_5{\partial}_{5,6},
\end{gather*}
be the linear parts of $\Pi_i^\lambda$ in a neighbourhood of
$(0,0,x_3,0,0,x_6)$. It follows that the Lie algebra determined by
$(\Pi_1^\lambda)^{(1)}$ (resp. $(\Pi_2^\lambda)^{(1)}$) is
isomorphic to $\rg_{3,\,-i\lambda}(\CC)$ (resp.
$\rg_{3,\,i\lambda}(\CC)$), considered as real Lie algebras.
\end{exam}

\medskip

\begin{rems}

\smallskip
\noindent $\ri$ Replacing $x_2$ by $-x_2$, the vector field $W_1$ is
transformed into $W_2$ while the other generators of the
distributions ${\mathcal C}_1^\lambda$ and ${\mathcal C}_2^\lambda$
remain unchanged. Thus, the corresponding leaves of both
distributions are diffeomorphic. However, these Lie-Poisson
structures are not equivalent since the associated Lie bialgebras
are non isomorphic. Furthermore, $\Pi_1^\lambda$ and
$\Pi_2^\lambda$, with $\lambda\neq 1$, are not equivalent as Poisson
structures since their linear parts give rise to non isomorphic Lie
algebras.

\smallskip

\noindent $\rii$ The Lie algebras $\rg_{3,\,-i\lambda}(\CC)$ and
$\rg_{3,\,i\lambda}(\CC)$, considered as six-dimensional real Lie
algebras, are real analitically degenerate. In fact, the Poisson
tensor
\[ \Pi=  \left(\Pi_1^\lambda\right)^{(1)} +(x_1 \partial _1 + x_2 \partial _2 )\wedge
(x_4 \partial _4 + x_5 \partial _5 )             \] is
non-linearizable since it has rank  $6$ almost everywhere, while its
linear part $\left(\Pi_1^\lambda\right)^{(1)}$ has rank at most $4$.
Therefore,  $\rg_{3,\,-i\lambda}(\CC)$ is degenerate. A similar
argument works for the Lie algebra $\rg_{3,\,i\lambda}(\CC)$.
\end{rems}

\

\subsection*{Acknowledgments} We would like to thank I. Dotti for
useful discussions and suggestions. We are also grateful to J.-P.
Dufour for providing the example exhibited in \S5.3, Remark (ii).
The first author would like to thank the hospitality at the
International Centre for Theoretical Physics (Trieste), where part
of the research for this paper was done. The authors were partially
supported by grants from ICTP (Italy), Conicet and Secyt U.N.C.
(Argentina).

\

\

\end{document}